\documentclass[reqno,11pt]{amsart}
\usepackage{latexsym}
\usepackage{amssymb,latexsym,amscd,amsmath,amsthm}
\usepackage{hyperref}
\def\R{\ensuremath{\mathbb R}}
\def\C{\ensuremath{\mathbb C}}
\def\Z{\ensuremath{\mathbb Z}}

\def\id{\operatorname{id}}
\let\eps=\varepsilon
\let\thet=\vartheta
\let\phy=\varphi
\def\rk{\operatorname{rk}}
\def\ind{\operatorname{ind}}

\def\rk{\operatorname{rk}}
\def\det{\operatorname{det}\nolimits}
\def\Ad{\operatorname{Ad}}
\def\ad{\operatorname{ad}}
\def\ch{\operatorname{ch}}
\def\aD{{\,\widetilde{\!\vphantom{x}\smash{\mathrm{ad}}\!}\,}}
\def\aDsl{{\,\widehat{\widetilde{\!\vphantom{x}\smash{\mathrm{ad}}\!}}\,}}

\def\adx_#1{\aD_{\frp,#1}}
\def\ady_#1{\aDsl_{\frp,#1}}
\def\pix{{\tilde\pi}}

\def\fre{{\mathfrak e}}
\def\frg{{\mathfrak g}}
\def\frh{{\mathfrak h}}

\def\frm{{\mathfrak m}}
\def\frp{{\mathfrak p}}
\def\frq{{\mathfrak q}}
\def\frs{{\mathfrak s}}
\def\frt{{\mathfrak t}}

\def\del{\partial}
\def\norm#1{\left\|#1\right\|}
\def\abs#1{\left|#1\right|}
\def\<{\langle}
\def\>{\rangle}
\def\Hom{{\operatorname{Hom}}}
\def\End{{\operatorname{End}}}
\def\Aut{{\operatorname{Aut}}}
\def\Vol{{\operatorname{vol}}}
\def\SU{\mathrm{SU}}
\def\U{\mathrm U}
\def\SO{\mathrm{SO}}
\def\Spin{\mathrm{Spin}}

\def\spin{\mathfrak{spin}}
\def\punkt{\,\cdot\,}

\def\Dsl{{\tilde D}}
\def\Adach{\hat A}
\def\Ldach{\hat L}
\def\Adachsl{\skew7\tilde{\skew6\hat A}}
\def\Ldachsl{\skew4\tilde{\skew3\hat L}}
\def\LC{^{\mathrm{LC}}}
\def\Cl{\operatorname{Cl}}
\def\sign{\operatorname{sign}}
\def\trace{\operatorname{tr}}
\def\strace{\operatorname{str}}
\def\even{{\mathrm{even}}}
\def\odd{{\mathrm{odd}}}
\numberwithin{equation}{section}
\swapnumbers
\theoremstyle{plain}
\newtheorem{Lemma}[equation]{Lemma}
\newtheorem{Proposition}[equation]{Proposition}
\newtheorem{Corollary}[equation]{Corollary}
\newtheorem*{thm}{Theorem}
\newtheorem{Theorem}[equation]{Theorem}

\theoremstyle{definition}

\theoremstyle{remark}

\newtheorem{Remark}[equation]{Remark}

\begin{document}

\title{Eta invariants of homogeneous spaces}
\author{S. Goette}
\address{Mathematisches Institut, Universit\"at Freiburg,
Eckerstra\ss e~1,
79104 Freiburg, Germany}
\email{sebastian.goette@math.uni-freiburg.de}
\begin{abstract}
We derive a formula for the $\eta$-invariants
of equivariant Dirac operators on quotients of compact Lie groups,
and for their infinitesimally equivariant extensions.
As an example, we give some computations for spheres.
\end{abstract}
\maketitle

Quotients~$M=G/H$ of compact Lie groups provide many important examples
of Riemannian manifolds with non-negative sectional curvature.
The primary characteristic classes and numbers of these spaces
have been computed by Borel and Hirzebruch in~\cite{BH}.

The $\eta$-invariant has been introduced by Atiyah, Patodi and Singer
in~\cite{APS}\ as a boundary contribution in an index theorem
for manifolds with boundary.
It can be used to construct certain secondary invariants
of compact manifolds~$M$
that were originally defined using zero-bordisms.
For example, the Eells-Kuiper and Kreck-Stolz invariants
distinguish homeomorphic homogeneous manifolds
that are not diffeomorphic.
These invariants can be expressed in terms of $\eta$-invariants
and Chern-Simons numbers, see e.g.~\cite{D1}, \cite{KS}.
Although the diffeomorphism type of many homogeneous manifolds~$G/H$
is well-known,
in some cases the explicit values of certain $\eta$-invariants are needed
to complete the diffeomorphism classification.
It is therefore worthwhile to have a formula for
the $\eta$-invariants of equivariant Dirac operators on homogeneous spaces.

First steps in this direction have been made
in~\cite{G1}, \cite{G2}, \cite{G3}.
There, we computed the equivariant $\eta$-invariant of
a different operator, called ``reductive'' or ``cubic'' Dirac operator,
and explained how to recover the $\eta$-invariant of the classical
Dirac operator.
However,
one complicated local term remained,
called ``Bott localisation defect'' below.
The central result of the present article is a formula for this defect term
that is similar to the formula for the equivariant $\eta$-invariant
of the reductive Dirac operator itself.
Thus we get a tractable formula for equivariant $\eta$-invariants
of homogeneous spaces
that is useful for explicit computations.
This formula has already been applied in a joint paper~\cite{GKS}\
with N. Kitchloo and K. Shankar
to calculate the Eells-Kuiper invariant
of the Berger space~$\SO(5)/\SO(3)$,
and to determine its diffeomorphism type.
Our formula can be summarised as follows;
details and notation will be explained later in the article.

\begin{thm}[see Theorem~\ref{MainTheorem}\ below]
  Let~$G\supset H$ be compact Lie groups,
  and let~$D^\kappa$ be the equivariant Dirac operator on~$M=G/H$
  twisted by the local bundle~$V^\kappa M$
  associated to an irreducible $\mathfrak h$-representation
  with highest weight~$\kappa$.
  Then~$\eta(D^\kappa)$ is a sum of the following terms:
  \begin{enumerate}
  \item a representation theoretic expression
    that vanishes if~$\rk G\ne\rk H+1$,
    that depends on the position of~$\mathfrak h$ in~$\mathfrak g$
    and on~$\kappa$,
    and that takes the form
    \begin{multline*}
	\kern3em
	2\sum_{w\in W_G}\frac{\sign(w)}{\delta(wX)}
	\,\Biggl(\Adach\bigl(\delta(wX)\bigr)
		\,e^{-\left(\alpha-\frac\delta2\right)(wX)}
		\,\prod_{\beta\in\Delta_G^+}\Adach\bigl(\beta(wX)\bigr)\\
	-e^{-(\kappa+\rho_H)(wX|_\frs)}
		\,\prod_{\beta\in\Delta_G^+}\Adach\bigl(\beta(wX|_\frs)\bigr)\Biggr)
		\,\prod_{\beta\in\Delta_G^+}\frac{-1}{\beta(X)}
		\Biggr|_{X=0}
    \end{multline*}
    if $\rk G=\rk H+1$,
  \item
    a local Chern-Simons theoretic contribution
	$$2\int_M\Adachsl\bigl(TM,\nabla^0,\nabla\LC\bigr)
			\,\ch\bigl(V^\kappa M,\nabla^\kappa\bigr)$$
    given by the integral of a constant multiple of the volume form of~$M$, and
  \item
    the integer
	$$\sum_{\gamma\in\hat G}\dim\bigl(V^\gamma\bigr)
		\,\Bigl(\eta\bigl({}^{\gamma\!} D^\kappa\bigr)
			-(\eta+h)\bigl({}^{\gamma\!}\Dsl^\kappa\bigr)\Bigr)$$
    arising as the spectral flow between two equivariant Dirac-type operators.
  \end{enumerate}
\end{thm}

We also give a reformulation of this result for the $\eta$-invariant
of the odd signature operator.
More generally,
we also obtain a formula for the infinitesimally equivariant
$\eta$-in\-var\-iant~$\eta_{\mathfrak g}(D^\kappa)$ of~\cite{G3};
this invariant is the universal $\eta$-form
for all families with fibrewise Dirac operator~$D^\kappa$
and compact structure group~$G$.
In fact,
we will compute the classical $\eta$-invariant~$\eta(D^\kappa)$
by evaluating~$\eta_{\mathfrak g}(D^\kappa)$ at~$X=0$.
Our main result is stated in Theorem~\ref{MainTheorem}\ for general equivariant
Dirac operators,
and in Corollary~\ref{BCor}\ for the slightly different special case
of the odd signature operator.
As an example,
we compute the infinitesimally equivariant $\eta$-invariants
of the untwisted Dirac operator and the odd signature operator
for round spheres,
and thus we obtain the corresponding $\eta$-forms of sphere bundles
with compact structure group.

Let us sketch our method.
In~\cite{G1}, \cite{G2}, we presented a formula
for the classical equivariant $\eta$-invariant of Slebarski's
deformed Dirac operator~$\Dsl^\kappa=D^{\frac13,\kappa}$
(\cite{Sl}, \cite{Ko}).
This invariant is a function~$\eta_G(\Dsl^\kappa)\colon G\to\C$
that is continuous
on the subset~$G_0\subset G$ of elements that act freely on~$M$,
and its value at the neutral element~$e\in G$ is just~$\eta(\Dsl^\kappa)$.
The singularity of~$\eta_G(\Dsl^\kappa)$ near~$e$ has implications
on fixpoint sets of $G$-manifolds~$N$ with~$\del N=M$,
see~\cite{G3}.
The difference between~$\eta_G(D^\kappa)$ and~$\eta_G(\Dsl^\kappa)$
is comparatively easy to control on~$G_0$,
but since~$\eta_G(D^\kappa)$ and~$\eta_G(\Dsl^\kappa)$ are not continuous
at~$e$,
we do not obtain the value of~$\eta(D^\kappa)$ by this method.

Instead,
we will use the infinitesimally equivariant
$\eta$-invariant~$\eta_{\mathfrak g}(D^\kappa)$ of~\cite{G3},
which is a power series on~$\frg$
with constant term~$\eta(D^\kappa)$.
For~$X\in\mathfrak g$ without zeros on~$M$,
the infinitesimal $\eta$-invariant~$\eta_X(D^\kappa)$
differs from the classical equivariant
$\eta$-invariant~$\eta_{e^{-X}}(D^\kappa)$
by the integral of a certain differential form on~$M$.
This way,
one obtains a formula for the power series~$\eta_{\mathfrak g}(D^\kappa)$,
and the classical $\eta$-invariant
is just its constant term.

However,
the integrand in the formula
for the difference~$\eta_X(D^\kappa)-\eta_{e^{-X}}(D^\kappa)$
stated in~\cite{G3}\ is in general not invariant under the action of~$G$,
so that one cannot reduce the problem
to a calculation at a single point in~$M$.
The new contribution in the present article is an evaluation of this
integral in terms of representation theoretic data
of the groups~$G$ and~$H$ and the relative position of their maximal tori,
up to an equivariant Chern-Simons term,
see Theorem~\ref{BLDThm}.
This is an improvement compared with respect to~\cite{G1},
\cite{G2}, \cite{G3}
for two reasons.
First, our formula for~$\eta(D^\kappa)$
involves no more equivariant differential forms.
Second,
for symmetric spaces,
the operators~$D^\kappa$ and~$\Dsl^\kappa$ are equal,
which means that there will be no Chern-Simons term and no spectral flow.

This paper is organised as follows.
In section~\ref{EtaKapitel},
we recall the results of~\cite{G1}, \cite{G2}, \cite{G3}
on equivariant $\eta$-invariants.
We also calculate the infinitesimally equivariant $\eta$-invariants of spheres.
In section~\ref{BottKapitel},
we give a formula for the Bott localisation defect on homogeneous spaces
and present our main result.

We wish to thank W. Soergel for some helpful comments.
We wish to thank an anonymous referee whose comments helped to make the paper
more readable.
Also,
we are indebted to K. Shankar and N. Kitchloo for their constant encouragement,
without which this paper would probably not have been written.

\section{Equivariant \texorpdfstring{$\eta$}{eta}-invariants}\label{EtaKapitel}

We recall some facts about $\eta$-invariants and homogeneous spaces
from~\cite{G1}, \cite{G2}, \cite{G3}.
In section~\ref{SymmAbschnitt},
we compute the infinitesimally equivariant $\eta$-invariant
for the untwisted Dirac operator and the signature operator on round spheres.

\subsection{Equivariant \texorpdfstring{$\eta$}{eta}-invariants
and their infinitesimal analogues}%
\label{DefAbschnitt}
Let~$(M,g^M)$ be an oriented Riemannian manifold,
and let~$G$ be a compact group that acts on~$M$ by isometries.
Let~$\mathcal E\to M$ be a $G$-equivariant Hermitian vector bundle,
equipped with a $G$-equivariant Clifford multiplication.
Let~$\nabla^{TM}$ denote a metric connection on~$M$,
and let~$\nabla^{\mathcal E}$ be a $G$-equivariant unitary connection on~$\mathcal E$
that satisfies the Leibniz rule
	$$\nabla^{\mathcal E}(v\cdot s)
	=\nabla^{TM} v\cdot s+v\cdot\nabla^{\mathcal E}s$$
for all vector fields~$v$ on~$M$ and all sections~$s\in\Gamma(\mathcal E)$.
Then~$\mathcal E$ is a {\em $G$-equivariant Clifford module\/} over~$M$,
and we have the  $G$-equivariant {\em Dirac operator\/}
$$\begin{CD}
	D\colon\Gamma(\mathcal E)
		@>\nabla^{\mathcal E}>>\Gamma(T^*M\otimes\mathcal E)
		@>g^M>>\Gamma(TM\otimes\mathcal E)
		@>\cdot>>\Gamma(\mathcal E)\;.
  \end{CD}$$
The most natural choice  for~$\nabla^{TM}$
is the Levi-Civita connection~$\nabla\LC$,
however,
on homogeneous spaces it is easier to work
with the connection~$\nabla^{\frac13}$,
see section~\ref{HomogenAbschnitt}.

Let~$g\in G$,
assume~$\Re s\ge s_0>0$,
and define
\begin{equation}\label{EqEtaFunction}
  \eta_g(D,s)
  =\sum_\lambda\sign(\lambda)\,\abs\lambda^{-s}
	\,\trace\bigl(g|_{E_\lambda}\bigr)
  =\int_0^\infty\frac{t^{\frac{s-1}2}}{\Gamma\bigl(s+\frac12\bigr)}
	\,\trace\bigl(gD\,e^{-tD^2}\bigr)\,dt\;,
\end{equation}
cf.~\cite{D2}\ and~\cite{Z1},
where~$E_\lambda$ denotes the Eigenspace of the Eigenvalue~$\lambda$.
The function~$\eta_g(D,s)$ admits a meromorphic continuation to~$\C$
that is finite at~$s=0$,
and the {\em equivariant $\eta$-invariant\/} of~$D$ at~$g$ is given by
\begin{equation}\label{EqEtaDef}
  \eta_g(D)=\eta_g(D,0)\;.
\end{equation}
If we have chosen~$\nabla^{TM}=\nabla\LC$,
then the second representation of~$\eta_G(D,s)$ actually converges
at~$s=0$.
The kernel of~$D$ does not contribute to~\eqref{EqEtaFunction};
we define
\begin{equation*}
  h_g(D)=\trace(g|_{\ker D})\;.
\end{equation*}

To define the infinitesimal analogue of~$\eta_g(D)$,
let~$X\in\frg$ be an element of the Lie algebra of~$G$,
and define the Killing field~$X_M$ on~$M$ by
\begin{equation}\label{XMDef}
  X_M|_p=\frac d{dt}\Bigr|_{t=0}e^{-tX}p\;.
\end{equation}
Let~$\mathcal L_X^{\mathcal E}$ denote the infinitesimal action of~$X$
on sections of~$\mathcal E$.
We take a Dirac operator~$D$ associated with the Levi-Civita connection,
and we construct deformations of~$D$ and~$D^2$ by
	$$D_X=D-\frac14\,X_M\cdot
		\qquad\text{and}\qquad
	H_X=D_{-X}^2+\mathcal L_X^{\mathcal E}\;,$$
where~$X_M$ acts by Clifford multiplication, see~\cite{BGV}.
Then the
{\em infinitesimal $\frg$-equivariant $\eta$-in\-var\-iant\/}~$\eta_\frg(D)
\in\C[\![\frg^*]\!]$
of~$D$ is defined as
\begin{equation}\label{InfEtaDef}
  \eta_X(D)
  =\int_0^\infty\frac1{\sqrt{\pi t}}
		\,\trace\Bigl(D_{\frac Xt}
			\,e^{-tH_{\frac Xt}}\Bigr)\,dt\;,
\end{equation}
see~\cite{G3}, cf.~\cite{BC1}, \cite{BC2}.
Note that it is conjecturally possible to define~$\eta_\frg(D)$
as a function of~$rX$ for~$X\in\frg$ and~$r\in\R$ sufficiently small,
rather than as a formal power series.
It is also possible to define a mixed equivariant
$\eta$-invariant~$\eta_{g,X}(D)$, where~$g\in G$ and~$X\in\frg$ commute;
this mixed invariant would then occur in an equivariant index theorem
for manifolds with boundary combining Theorem~\ref{APSDThm}~(1) and~(2),
and in an equivariant index theorem for families with compact structure group,
cf.\ Remark~\ref{APSDBCRem}.
However,
this extra generality is not required for the applications we have in mind.

The equivariant $\eta$-invariant and its infinitesimal cousin
are key ingredients in theorem~\ref{APSDThm} below,
which is due to Atiyah, Patodi, Singer~\cite{APS}\ and Donnelly~\cite{D2},
and its infinitesimal analogue in~\cite{G3},
which is related to the family index theorem of Bismut and Cheeger
in~\cite{BC2}.
Suppose that~$M$ is the boundary of a compact oriented manifold~$N$,
which has a collar isometric to~$M\times[0,\eps]$.
Suppose that there exists a Clifford module~$\overline{\mathcal E}
=\smash{\overline{\mathcal E}}^+\oplus\smash{\overline{\mathcal E}}^-\to N$,
such that~$\smash{\overline{\mathcal E}}^+|_M\cong\mathcal E$
with compatible Clifford multiplications.
Then the Dirac operator~$D$ is closely related
to the Dirac operator~$\overline D$ acting on sections of~$\overline{\mathcal E}$.

We also need equivariant characteristic differential forms.
The {\em equivariant Riemannian curvature\/} is defined as
	$$R_X
	=\bigl(\nabla\LC-2\pi i\,\iota_{X_M}\bigr)^2+2\pi i\,\mathcal L_X^{TM}
	=R+2\pi i\,\mu\LC_X\;,$$
where~$\iota_{X_M}$ denotes interior multiplication of a differential form
by~$X_M$,
and~$\mu\LC_X=\mathcal L_X^{TM}-\nabla\LC_{X_M}$ is the {\em Riemannian moment.\/}
Similarly,
define
	$$F^{\mathcal E}_X
	=\bigl(\nabla^{\mathcal E}-2\pi i\,\iota_{X_M}\bigr)^2+2\pi i\,\mathcal L_X^{\mathcal E}
	=F^{\mathcal E}+2\pi i\,\mu^{\mathcal E}_X\;.$$
Following~\cite{BGV},
the {\em equivariant twisting curvature\/} of~$\mathcal E$ is defined as
\begin{equation}\label{FESDef}
  F^{\mathcal E/\mathcal S}_X
  =F^{\mathcal E}_X-\frac14\sum_{i,j}\<R_Xe_i,e_j\>
			\,e_i\cdot e_j\cdot\mathord{}\;,
\end{equation}
where~$e_1$, \dots, $e_n$ is a local orthonormal frame of~$TM$,
which acts on~$\mathcal E$ by Clifford multiplication.
Then we have the equivariant characteristic differential forms
\begin{align*}
  \Adach_X\bigl(TM,\nabla\LC\bigr)
  &=\det^{\textstyle\frac12}
	\biggl(\frac{R_X/4\pi i}{\sinh(R_X/4\pi i)}\biggr)\\
	\text{and}\qquad
  \ch_X\bigl(\mathcal E/\mathcal S,\nabla^{\mathcal E}\bigr)
  &=\trace\biggl(e^{-\textstyle\frac{F^{\mathcal E/\mathcal S}_X}{2\pi i}}
	\biggr)\;.
\end{align*}
Recall that the equivariant exterior derivative is defined as
	$$d_X=d-2\pi i\,\iota_{X_M}\;,$$
such that~$d_X^2+2\pi i\,\mathcal L_X=0$.
Then the forms~$\Adach_X(TM,\nabla\LC)$
and~$\ch_X(\mathcal E/\mathcal S,\nabla^{\mathcal E})$ are
$G$-invariant, equivariantly closed,
and independent of the choice of an equivariant connection
up to equivariantly exact forms.
Note that in contrast to~\cite{BGV} and~\cite{G3},
we use the classical convention of including powers
of~$2\pi i$ in all definitions above.
Let~$N_g$ denote the fixpoint set of~$g$ on~$N$.

Let~$\ind_g(\overline D)=\strace(g|_{\ker D})$ denote the equivariant index
of~$\overline D$ (regarded as the character of a virtual $G$-representation
at~$g\in G$)
with respect to the APS boundary conditions.

\begin{Theorem}[\cite{APS}, \cite{D2}, \cite{G3}]\label{APSDThm}
  Let~$\alpha_g$ denote the characteristic form on~$N_g$
  given by the constant term in the asymptotic development
  of~$\strace\bigl(g\,e^{-t\overline D^2}\bigr)$ for~$t\to 0$,
  then
  \begin{align*}
    \ind_g\bigl(\overline D\bigr)
    &=\int_{N_g}\alpha_g(D)
	-\frac{\eta_g(D)+h_g(D)}{2}\;.\tag1
  \end{align*}
  Assume moreover that~$\overline D$ is associated
  to the Levi-Civita connection,
  then so is~$D$,
  and
  \begin{align*}
    \ind_{e^{-X}}\bigl(\overline D\bigr)
    &=\int_N\Adach_X\bigl(TM,\nabla\LC\bigr)
	\,\ch_X\bigl(\mathcal E/\mathcal S,\nabla^{\mathcal E}\bigr)
	-\frac{\eta_X(D)+h_{e^{-X}}(D)}2\;.\tag2
  \end{align*}
\end{Theorem}

\begin{Remark}\label{APSDBCRem}
  Theorem~\ref{APSDThm}~(2) is in fact a special case
  of the Bismut-Cheeger index theorem for families of manifolds with boundary
  in~\cite{BC2}.
  Suppose that~$P\to B$ is a $G$-principal bundle
  with curvature form~$\omega$ and curvature~$\Omega$.
  If we consider the family~$P\times_GN\to B$
  with boundary~$P\times_GM\to B$,
  then Bismut-Cheeger's theorem is equivalent to~(2) above
  by the Chern-Weil construction.
  In particular,
  we can recover the $\eta$-form on~$B$ of the family of Dirac operators
  induced by~$D$ on~$M$ from~$\eta_X(D)$.
\end{Remark}

We now recall the relation between the two equivariant $\eta$-invariants
introduced in~\eqref{EqEtaDef} and~\eqref{InfEtaDef} above.
From the dual of the Killing field~$X_M$,
we construct an equivariant differential form
\begin{equation}\label{ThetaXFormel}
  \thet_X=\frac1{2\pi i}\,\<X_M,\punkt\>\;.
\end{equation}
For~$X\ne 0$,
we have an $L_1$-current
	$$\frac{\thet_X}{d_X\thet_X}
	=-\int_0^\infty{\thet_X}
		\,e^{t\,{d_X\thet_X}}\,dt$$
on~$M$ by Proposition~\ref{LeinsProp}\ below.
Note that~$\frac{\thet_X}{d_X\thet_X}$ has a singularity
at the zero set of~$X$,
but for a smooth compactly supported form~$\alpha$ on~$M$,
the integral~$\int_M\frac{\thet_X}{d_X\thet_X}\alpha$ is still well-defined,
hence the term ``current''.
The current~$\frac{\thet_X}{d_X\thet_X}$ will be referred to as the
{\em Bott localisation defect\/} because of~\eqref{BVLocFormel} below.

\begin{Theorem}[\cite{G3}, Theorem~0.5]\label{GdreiThm}
If the Killing field~$X_M$ has no zeros on~$M$,
then
	$$\eta_{rX}(D)
	=\eta_{e^{-rX}}(D)+
		2\int_M\frac{\thet_{rX}}{d_{rX}\thet_{rX}}
			\,\Adach_{rX}(M)\,\ch_{rX}(\mathcal E/\mathcal S)
		\quad\in\C[\![r]\!]\;.$$
\end{Theorem}


If there exists a $G$-manifold manifold~$N$ that bounds~$M$ equivariantly
and if there exists a Dirac operator~$\overline D$ on~$N$ related to~$D$
as above,
then~Theorem~\ref{GdreiThm}\ can be proved by comparing~Theorem~\ref{APSDThm}~(1) and~(2) using
Kalkman's localisation formula for manifolds with boundaries~(\cite{K}).
A general argument is given in~\cite{G3}.

We will give a formula for the Bott localisation defect
in section~\ref{BottKapitel}
in the special case that~$M$ is a homogeneous space of compact type.

\subsection{Homogeneous spaces}\label{HomogenAbschnitt}
We recall the formula for the equivariant $\eta$-invariants of reductive
Dirac operators,
and its relation to the $\eta$-invariants of Riemannian Dirac operators.

Let~$M=G/H$ be a quotient of compact Lie groups
with Lie algebras~$\frh\subset\frg$.
Recall that all $G$-equivariant vector bundles
over~$M=G/H$
are of the form
	$$V^\kappa M=G\times_\kappa V^\kappa\to M\;,$$
where~$(\kappa,V^\kappa)$ is a representation of~$H$.
We identify sections~$s$ of~$V^\kappa M$ with $H$-equivariant functions
\begin{equation*}
  \hat s\colon G\to V^\kappa
	\qquad\text{with}\qquad
  \hat s(gh)=\kappa_h^{-1}\,\hat s(g)
\end{equation*}
such that~$s(gH)=[g,\hat s(g)]$.

We fix an $\Ad_G$-invariant metric on~$\frg$
and let~$\frm=\frh^\perp\subset\frg$.
Let~$\pi=\Ad|_{H\times\frm}$ denote the isotropy representation
of~$H$ on~$\frm$.
Then the tangent bundle~$TM$ is isomorphic to~$G\times_\pi\frm$
via
\begin{equation}\label{TMiso}
  [g,v]=\frac d{dt}ge^{tv}
\end{equation}
and carries an induced
{\em normally homogeneous\/} Riemannian metric.

On most bundles~$V^\kappa M$,
we will mainly use the reductive connection~$\nabla^0=\nabla^{0,\kappa}$,
which can be written as
\begin{equation}\label{NablaNullDef}
	\widehat{\nabla^0_Vs}=\hat V(\hat s)\;.
\end{equation}
If~$\kappa$ is a unitary (orthogonal) representation,
then~$\nabla^{0,\kappa}$ is a unitary (orthogonal) connection.
A subscript~$\frh$ or~$\frm$ will denote orthogonal projection
to that subspace of~$\frg$.
One easily calculates the curvature of~$\nabla^{0,\kappa}$ as
\begin{equation}\label{CurvatureNullFormel}
  \widehat{F^{0,\kappa}_{V,W}s}
  =-\kappa_{*\left[\hat V,\hat W\right]_\frh}\hat s\;.
\end{equation}

We consider a family of connections on~$TM$
given by
\begin{equation}\label{NablaTDef}
  \widehat{\nabla^t_VW}
  =\hat V(\hat W)+t\,\bigl[\hat V,\hat W\bigr]_\frm\;.	
\end{equation}
For~$t=0$,
we obtain again the reductive connection.
Note that because the Lie bracket of vector fields on~$M$ is given by
	$$\widehat{[V,W]}
  	=\hat V\bigl(\hat W\bigr)-\hat W\bigl(\hat V\bigr)
		+\bigl[\hat V,\hat W\bigr]_\frm\;,$$
the Levi-Civita connection on~$TM$ with respect to a normal metric
is~$\nabla\LC=\nabla^{0,\pi}+\frac12\,[\punkt,\punkt]_\frm$.

An equivariant Clifford module over~$M$
is a $G$-equivariant vector bundle~$\mathcal E\to M$
together with a $G$-equivariant
Clifford multiplication~$TM\times\mathcal E\to\mathcal E$.
We can regard the isotropy representation
as a homomorphism~$\pi\colon\frh\to\spin(\frm)$,
where~$\spin(\frm)$ denotes the spin group
associated to the vector space~$\frm$.
Let~$n$ be the dimension of~$M$, and
let~$\pix$ denote the pullback of the spin representation of~$\spin(\frm)$
to~$\frh$,
which acts on a complex $2^{\left[\frac n2\right]}$-dimensional
vector space~$S$.
If~$\pix$ integrates to an $H$-representation,
then~$\mathcal S=G\times_HS$ is the $G$-equivariant spinor bundle
(unique if~$H$ is connected),
which is then the most elementary Clifford bundle.
In general,
every equivariant Clifford module is of the form
\begin{equation}\label{SkappaDef}
  \mathcal S^\kappa M=G\times_H\bigl(S\otimes W^\kappa\bigr)\to M\;,
\end{equation}
where~$\kappa$ is an $\frh$-representation such that~$\pix\otimes\kappa$
integrates to an $H$-representation (\cite{G2}, Lemma~3.4).
For example,
if~$M$ is even-dimensional,
the complexified bundle of exterior differential forms on~$M$
is precisely the $G$-equivariant Clifford module
	$$\Lambda^*TM\otimes\C
	=\mathcal S^\pix M
	=G\times_H(S\otimes S)\to M\;.$$

Let~$\mathcal S^\kappa M\to M$ be a $G$-equivariant Clifford module.
Following Slebarski~\cite{Sl},
we define a family of Dirac operators on~$\Gamma(\mathcal S^\kappa M)$.
For~$X\in\frg$,
we define~$\adx_X\in\End S$ by
	$$\adx_X
	=\frac14\,\sum_{i,j}\<[X,e_i],e_j\>\,e_i\cdot e_j\cdot\mathord{}\;,$$
where~$e_i$, $e_j$ run through an orthonormal base of~$\frm$.
Note that~$\adx_X=\pix_{*X}$ for~$X\in\frh$.
Then we set
	$$\widehat{D^{t,\kappa}s}
	=\sum_ie_i\cdot\bigl(e_i(\hat s)+t\,\adx_{e_i}\bigr)\;.$$
For~$t=\frac12$,
the operator~$D^\kappa=D^{\frac12,\kappa}$
is associated to the Levi-Civita connection on~$\mathcal S M$
and the reductive connection on~$V^\kappa M$
(but these two bundles are in general only locally well defined).

The operator~$\Dsl^\kappa=D^{\frac13,\kappa}$ has distinguished properties.
It was shown in~\cite{Sl} and~\cite{G1}, \cite{G2} that~$\Dsl^\kappa$
exhibits a very natural behaviour with respect to homogeneous fibrations,
hence it was called the ``reductive'' Dirac operator in~\cite{G2}
(although it is not constructed from the reductive connection on~$S^\kappa M$).

In order to state the formula for~$\eta_G(\Dsl^\kappa)$,
we need some more notation and conventions.
We choose maximal tori~$S\subset T$ of~$H\subset G$
with Lie algebras~$\frs\subset\frt$,
and fix Weyl chambers~$P_G\subset i\frt^*$ and~$P_H\subset i\frs^*$.
Let~$\Delta_G^+$ and~$\Delta_H^+$ denote
the corresponding sets of positive roots,
and let~$\rho_G$ and~$\rho_H$ be their half sums.
The choices of~$P_G$ and~$P_H$ also determine orientations on~$\frg/\frt$
and~$\frh/\frs$ as follows.
If~$\beta_1$, \dots, $\beta_l\in i\frt^*$ are the positive roots of~$\frg$
with respect to~$P_G$,
we can choose a complex structure on~$\frg/\frt$
and a complex basis~$z_1$, \dots, $z_n$
such that~$\ad|_{\frt\times(\frg/\frt)}$ takes the form
\begin{equation}\label{GTOrientation}
  \ad_X=\begin{pmatrix}\beta_1(X)\\&\ddots\\&&\beta_l(X)\end{pmatrix}
		\qquad\text{for all~$X\in\frt$.}
\end{equation}
Then we declare the basis~$z_1$, $i\,z_1$, $z_2$, \dots, $i\,z_l$
of~$\frg/\frt$ as a real vector space
to be positive oriented.
The orientation of~$H/S$ is constructed similarly.

If we fix an orientation on~$\frm=\frg/\frh$
and choose orientations on~$\frg/\frt$ and~$\frh/\frs$ as above,
there is a unique orientation on~$\frt/\frs$
such that the two induced orientations on
\begin{equation}\label{OrientFormel}
  \frg/\frs\cong\frm\oplus(\frh/\frs)\cong(\frg/\frt)\oplus(\frt/\frs)
\end{equation}
agree.
We assume for the moment that~$\rk G=\rk H+1$,
so~$\frs^\perp\subset\frt$ is one-dimensional.
Let~$E\in\frt/\frs\cong\frs^\perp$ be the positive unit vector,
and let~$\delta\in i\frt^*$ be the unique weight such that
\begin{equation}\label{DeltaDef}
  -i\delta(E)>0\qquad\text{and}\qquad\delta(X)\in2\pi i\,\Z\iff e^X\in S
\end{equation}
for all~$X\in\frt$.

By abuse of notation,
let~$\kappa\in i\frs^*$ denote the highest weight
of the $\frh$-representation~$\kappa$ used to construct~$\mathcal S^\kappa M$.
Then there is a unique weight~$\alpha\in i\frt^*$ of~$\frg$ such that
\begin{equation}\label{AlphaDef}
  \alpha|_\frs=\kappa+\rho_H
	\qquad\text{and}\qquad
  -i(\alpha-\delta)(E)<0\le-i\alpha(E)\;.
\end{equation}
Let~$A_G$ denote the alternating sum over the Weyl group of~$G$
acting on~$\frt$,
and write~$A_G(\rho_G)$ shorthand for~$A_G(e^{\rho_G(\punkt)})$.
Recall that for~$X\in\frt$, $e^X\in T$ is regular iff~$A_G(\rho_G)(X)\ne0$.
We can now compute the classical equivariant eta-invariant of~$\Dsl^\kappa$.

\begin{Theorem}[\cite{G1}, \cite{G2}]\label{EtaThm}
  If~$\rk G\ne\rk H+1$,
  then~$\eta_G(\Dsl^\kappa)=0$ identically .
  If~$\rk G=\rk H+1$,
  then~$\eta_G(\Dsl^\kappa)$ is continuous on the set~$T_0\subset T$
  that acts freely on~$M$.
  Moreover,
  if~$\alpha$, $\delta\in i\frs^*$ and~$E\in\frt$ are given as above,
  then for all~$X\in\frt$ such that~$e^X\in T_0$ and~$e^X$ is regular,
	$$(\eta+h)_{e^X}\bigl(\Dsl^\kappa_M\bigr)
	=\frac{A_G\Bigl(\frac 1{\sinh\left(\frac\delta2(\punkt)\right)}
			\,e^{(\alpha-\frac\delta2)(\punkt)}\Bigr)}
		{A_G(\rho_G)}(X)\;.$$
\end{Theorem}

\subsection{Round Spheres}\label{SymmAbschnitt}
As an example for our main result,
we compute the infinitesimally equivariant $\eta$-invariant
of the untwisted Dirac operator and the odd signature operator
for odd-dimensional spheres.
Note that for a symmetric space,
the last two terms in Theorem~\ref{MainTheorem}~(1) and~(2) vanish,
so in particular,
we do not have to integrate equivariant Chern-Simons classes.

Fix~$n$ and embed~$G=\Spin(2n)\subset\Cl(2n)$.
We let~$H=\Spin(2n-1)$
denote the subgroup belonging to the Clifford subalgebra spanned by the
vectors~$e_1$, \dots, $e_{2n-1}$ 
of an orthonormal base of~$\R^{2n}$.
Then we can write~$S^{2n-1}=G/H$.

We choose the maximal tori~$S\subset H$ and~$T\subset G$
with Lie algebras
\begin{align*}
  \frt
  &=\textstyle
	\bigl\{\,X=\frac{x_1}2\,e_1e_2+\dots
		+\frac{x_n}2\,e_{2n-1}e_{2n}
	\bigm|(x_1,\dots x_n)\in\R^n\bigr\}\\
  \text{and}\qquad
  \frs
  &=\{\,X\in\frt\mid x_n=0\,\}\;.
\end{align*}
We also choose Weyl chambers
\begin{align*}
  P_G
  &=\bigl\{\,i\,(\gamma_1x^1+\dots+\gamma_nx^n)
	\bigm|\gamma_1\ge\dots\ge\gamma_{n-1}
		\ge\abs{\gamma_n}\,\bigr\}\\
  \text{and}\qquad
  P_H
  &=\bigl\{\,i\,(\gamma_1x^1+\dots+\gamma_{n-1}x^{n-1})
	\bigm|\gamma_1\ge\dots\ge\gamma_{n-1}\ge 0\,\bigr\}\;.
\end{align*}
Then we have
\begin{align}\begin{split}\label{WurzelFormel}
    \Delta_G^+
    &=\bigl\{\,i\,(x_j\pm x_k)\bigm|1\le j<k\le n\,\bigr\}\\
    \text{and}\qquad
    \Delta_H^+
    &=\bigl\{\,i(x_j\pm x_k)\bigm|1\le j<k\le n-1\bigr\}
	\cup\bigl\{\,ix_j\bigm|1\le j\le n-1\,\bigr\}\;.
\end{split}\end{align}
We fix~$\delta(X)=ix_n$ in~\eqref{DeltaDef},
so that~$S^{2n-1}$ is now oriented by~\eqref{GTOrientation}
and~\eqref{OrientFormel}.

\begin{Theorem}\label{OGThm}
  For~$X=(x_1,\dots,x_n)\in\frt$ as above,
  the infinitesimally equivariant $\eta$-invariant
  of the untwisted Dirac operator~$D$
  and the odd signature operator~$B=D^{\pix}$
  of the odd spheres are given by
  \begin{align*}
    \eta_X\bigl(D_{S^{2n-1}}\bigr)
    &=\frac{i^n}{2^{n-1}\sin\frac{x_1}2\cdots\sin\frac{x_n}2}
	\biggl(1-\sum_{j=1}^n\frac2{x_j}\,\sin\frac{x_j}2
		\prod_{k\ne j}\frac{x_k^2}{x_k^2-x_j^2}\biggr)\;,\tag1\\
    \eta_X\bigl(B_{S^{2n-1}}\bigr)
    &=i^n\,\cot\frac{x_1}2\cdots\cot\frac{x_n}2
	\biggl(1-\sum_{j=1}^n\frac2{x_j}\,\tan\frac{x_j}2
		\prod_{k\ne j}\frac{x_k^2}{x_k^2-x_j^2}\biggr)\;.\tag2
  \end{align*}
\end{Theorem}

Note that these formulas have been proved by Zhang for~$S^1$ in~\cite{Z2}
and for~$S^3$ by the author in~\cite{G3}.
For other symmetric spaces,
one can similarly compute~$\eta_X(D)$ and~$\eta_X(B)$
using the formulas for the equivariant $\eta$-invariants
in~\cite{G1}, \cite{G2}.

\begin{proof}
  The classical equivariant $\eta$-invariants of~$B$ and~$D$
  have been calculated in~\cite{APS}\ and~\cite{HR}.
  Let us start with the untwisted Dirac operator.
  We take~$X\in\frt$ as above.
  By~\cite{HR},
  we have
	$$\eta_{e^{-X}}\bigl(D_{S^{2n-1}}\bigr)
	=\frac{i^n}{2^{n-1}\,\sin\frac{x_1}2\cdots\sin\frac{x_n}2}\;.$$
  Let~$W_\pi^+$ denote the positive weights
  of the isotropy representation~$\pi$,
  which are all of multiplicity one in our case.
  With~\eqref{WurzelFormel} and~$\delta(X)=ix_n$,
  we calculate
	$$\frac1{\delta(-X)}
		\,\prod_{\beta\in\Delta_G^+}
			\frac{\beta(-X|_\frs)}{\beta(-X)}
		\,\prod_{\beta\in W_\pi^+}
			\frac1{2\sinh\bigl(\frac\beta2(-X|_\frs)\bigr)}
	=\frac i{x_n}\prod_{k=1}^{n-1}\frac{x_k^2}{(x_k^2-x_n^2)}
		\prod_{k=1}^{n-1}\frac i{2\sin\frac{x_k}2}\;.$$
  Note that the sum~$S_{G/H}$ over~$W_G/W_H$ of~\eqref{SGHDef}
  contains~$n$ similar terms,
  where~$x_n$ above is replaced by~$x_j$ for~$j=1$, \dots, $n$.
  Because the last two terms in Theorem~\ref{MainTheorem}~(2) drop out
  for a symmetric space,
  we get~(1).
  The proof of~(2) is similarly based on the formula for~$\eta_{e^X}(B)$
  in~\cite{APS} and Corollary~\ref{BCor}.
\end{proof}

\section{The Bott localisation defect}\label{BottKapitel}

The difference between
the classical equivariant $\eta$-invariant of~\cite{D2}
and the infinitesimal equivariant $\eta$-invariant of~\cite{G3}\
can be expressed in terms of the Bott localisation defect
as in~Theorem~\ref{GdreiThm}, see~\cite{G3}.
In this section,
we prove a fibration formula for the Bott localisation defect,
which is then applied to homogeneous spaces.

\subsection{A fibration formula}\label{CurrentKapitel}
Let~$\mathcal N_{M_X,M}$ denote the normal bundle to the zero set~$M_X$
of the Killing field~$X_M$ on~$M$,
and let~$e_X\bigl(\mathcal N_{M_X,M}\bigr)$ denote its equivariant Euler class.
Bott's localisation formula in equivariant cohomology
in the local version of~\cite{BGV}\ states that any equivariantly closed
form~$\alpha_X\in\Omega_\frg(M)$ is equivariantly cohomologous to the current
	$$\frac{\delta_{M_X}\,\alpha_X}{e_X\bigl(\mathcal N_{M_X,M}\bigr)}$$
for small~$X\in\frg$.
Here, $\delta_{M_X}$ denotes the $\delta$-distribution
at the fixpoint set of~$X$.
It follows from Berline-Vergne's proof in~\cite{BGV}\ that in fact
\begin{equation}\label{BVLocFormel}
  \alpha_X-\frac{\delta_{M_X}\,\alpha_X}{e_X\bigl(\mathcal N_{M_X,M}\bigr)}
  =d_X\biggl(\frac{\thet_X}{d_X\thet_X}\,\alpha_X\biggr)\;,
\end{equation}
where~$\thet_X=\frac1{2\pi i}\,X_M^*$ as in~\eqref{ThetaXFormel}.
The current~$\frac{\thet_X}{d_X\thet_X}$ also appears
in a localisation formula for manifolds with boundary in~\cite{K}.
In this section,
we state some properties of~$\frac{\thet_X}{d_X\thet_X}$.
In particular,
we give a formula for fibre bundles.
Let us begin with the following observation,
which holds for arbitrary Riemannian manifolds~$M$ with an isometric action
by a Lie group~$G$.

\begin{Proposition}\label{LeinsProp}
  If we set~$\frac{\thet_X}{d_X\thet_X}|_{M_X}=0$,
  then the current~$\frac{\thet_X}{d_X\thet_X}$ is locally of class~$L_1$.
  In particular,
  if~$\alpha\in\Omega^*(M)$ is continuous and compactly supported,
  then
	$$\int_M\frac{\thet_X}{d_X\thet_X}\,\alpha
	=-\int_0^\infty
		\biggl(\int_M\thet_X\,e^{t\,d_X\thet_X}\,\alpha\biggr)\,dt\;.$$
  In particular, both integrals exist.
\end{Proposition}

\begin{proof}
Because
	$$d_X\thet_X
	=\frac1{2\pi i}\,dX_M^*-\norm{X_M}^2\;,$$
it is clear that~$\frac{\thet_X}{d_X\thet_X}$ is integrable
on each compact subset of~$M\setminus M_X$.
Also,
the equation above holds for all~$\alpha$
with compact support in~$M\setminus M_X$.

Because~$X_M$ is a Killing field,
the zero set~$M_X$ is a totally geodesic submanifold of~$M$.
Let~$\mathcal N\to M_X$ denote the normal bundle,
then the Levi-Civita connection on~$M$
restricts to a connection~$\nabla^{\mathcal N}$ on~$\mathcal N$,
and the infinitesimal action~$\mu_X^{\mathcal N}$ of~$X$
is parallel with respect to~$\nabla^{\mathcal N}$.
The normal exponential map~$\exp^\perp\colon\mathcal N\to M$
is a local diffeomorphism near the zero section in~$\mathcal N$,
and
	$$e^{tX}\exp^\perp V=\exp^\perp\bigl(e^{t\mu_X^{\mathcal N}}V\bigr)$$
for all vectors~$V\in\mathcal N$.

Let~$\mathcal R$ denote the vector field near~$M_X$ given by
	$$\mathcal R|_{\exp^\perp V}=\frac d{dt}\Bigr|_{t=1}\exp^\perp tV\;,$$
and set~$r=\abs{\mathcal R}$ near~$M_X$.
To estimate the behaviour of~$dX_M^*$ and~$\norm{X_M}^2$ near~$M_X$,
we need some facts.

Because~$X_M$ is a Killing field,
it satisfies $\norm{X_M}^2>cr^2$ near~$M_X$ for some~$c>0$.
Cartan's formula for the exterior derivative implies that
	$$\bigl(dX_M^*\bigr)(V,W)=2\,\<\nabla_VX,W\>\;.$$
Now let~$V$ be a vector field on~$M$ near~$M_X$
that is tangential to~$M_X$
and parallel along all radial geodesics emanating from
and perpendicular to~$M_X$.
Let~$Y$, $Z$ be a vector fields normal to~$M_X$ with~$\nabla_YZ=0$,
then
\begin{equation}\label{ThetaXEstimates}
  \begin{gathered}
    \<X_M,V\>|_{M_X}=0\;,\qquad
    \bigl(Y\<X_M,V\>\bigr)|_{M_X}=\<\mu_X^{\mathcal N}V,Y\>=0\;,\\
	\text{and}\qquad
    Y\bigl(Z\<X_M,V\>\bigr)|_{M_X}
    =\bigl\<R_{Y,V}X_M-\bigl(\nabla_V\mu_X^{\mathcal N}\bigr)(Y),
	Z\bigr\>|_{M_X}
    =0
  \end{gathered}
\end{equation}
implies that
	$$\iota_V\,dX_M^*\le Cr^2\qquad\text{and}\qquad\<X_M,V\>\le Cr^3$$
for some~$C>0$.

We choose a local orthonormal frame~$e_1$, \dots, $e_n$ near~$M_X$
that is parallel along radial geodesics emanating from~$M_X$,
such that~$e_1$, \dots, $e_{2k}$ are normal to~$M_X$
and~$e_{2k+1}$, \dots, $e_n$ are tangential.
Then
\begin{equation}\label{ThetaXDThetaXEst}
  \frac{\thet_X}{d_X\thet_X}
  =-\frac1{2\pi i}\,\frac{X_M^*}{\norm{X_M}^2-\frac{dX_M^*}{2\pi i}}
  =-\sum_{j=0}^{\left[\frac{n-1}2\right]}(2\pi i)^{-j-1}
	\,\frac{X_M^*\,(dX_M^*)^j}{\norm{X_M}^{2j+2}}
  =O(r^{1-2k})
\end{equation}
near~$M_X$,
so the left hand side integral in Proposition~\ref{LeinsProp}\ exists.

It also follows from~\eqref{ThetaXDThetaXEst} that we may write
\begin{equation}\label{LeinsIntegralFormel}
  -\int_0^\infty\biggl(\int_M\thet_X\,e^{t\,d_X\thet_X}\,\alpha\biggr)\,dt
  =\int_M\frac{\thet_X}{d_X\thet_X}\,\alpha
	-\lim_{t\to\infty}\int_M\frac{\thet_X}{d_X\thet_X}
	\,e^{td_X\thet_X}\,\alpha\;.
\end{equation}
Let~$U_t$ be the tubular neighbourhood of~$M_X$
of radius~$t^{-\frac14}$ for~$t\gg 0$.
Because~$\norm{X_M}^2\ge cr^2$,
one has
\begin{equation}\label{LeinsZweiEst}
  \lim_{t\to\infty}\int_{M\setminus U_t}\frac{\thet_X}{d_X\thet_X}
	\,e^{td_X\thet_X}\,\alpha
  =0\;,
\end{equation}
so the rightmost term in~\eqref{LeinsIntegralFormel} localises near~$M_X$.

We identify~$U_t$
with the neighbourhood~$\mathcal N_t$ of radius~$t^{\frac14}$
by a rescaled normal exponential map
\begin{equation}\label{PhyTDef}
  v\mapsto\phy_t(v)=\exp^\perp\bigl(t^{-\frac1{ 2}}v\bigr)\;.
\end{equation}
Using~\eqref{ThetaXEstimates},
one finds that if~$A=\mu^{\mathcal N}_X\in\Gamma(\End\mathcal N)$,
then as~$t\to\infty$,
\begin{equation}\label{ThetaXLimit}
  \begin{gathered}
    t\,\phy_t^*\,X_M^*\to\<A\mathcal R,\punkt\>\;,\qquad
    t\,\phy_t^*\,dX_M^*\to \<2A\punkt,\punkt\>\;,\\
	\text{and}\qquad
    t\,\phy_t^*\norm{X_M}^2\to\norm{A\mathcal R}^2
  \end{gathered}
\end{equation}
uniformly on~$\mathcal N_t$,
with respect to the canonical metric on the total space of~$\mathcal N\to M_X$.
Moreover,
	$$\lim_{t\to\infty}\phy_t^*\alpha=\pi^*\alpha$$
is the pullback of~$\alpha$ by the bundle projection~$\pi\colon\mathcal N\to M_X$,
because~$\alpha$ is of class~$C^0$.
Using~\eqref{LeinsZweiEst}, this allows us to compute
\begin{multline}\label{LeinsLimit}
  \lim_{t\to\infty}\int_M\frac{\thet_X}{d_X\thet_X}
	\,e^{td_X\thet_X}\,\alpha\\
  \begin{aligned}
    &=-\lim_{t\to\infty}\int_{\mathcal N_t}\sum_{j,k}(2\pi i)^{-j-k-1}
	\,\frac{t\phy_t^*X_M^*\,(t\,\phy_t^*dX_M^*)^{j+k}}
		{\bigl(t\,\phy_t^*\norm{X_M}^2\bigr)^{j+1}\,k!}
	\,e^{-t\,\phy_t^*\norm{X_M}^2}
	\,\phy_t^*\alpha\\
    &=-\lim_{t\to\infty}
	\int_{\mathcal N}\sum_{j,k}(2\pi i)^{-j-k-1}
	  \,\frac{\<A\mathcal R,\punkt\>\,\<2A\punkt,\punkt\>^{j+k}}
		{\norm{A\mathcal R}^{2j+2}\,k!}
	  \,e^{-\norm{A\mathcal R}^2}\,\pi^*\alpha\\
    &=0\;.
  \end{aligned}
\end{multline}
This proves the existence of the right hand side integral
in Proposition~\ref{LeinsProp}.
Together with~\eqref{LeinsIntegralFormel},
we obtain the claimed equality.
\end{proof}

\begin{Theorem}\label{FibreCurrentTheorem}
  Let~$p\colon M\to B$ be a proper Riemannian submersion,
  and let~$G$ be a compact Lie group of isometries of~$M$
  that map fibres to fibres.
  Then
	$$\frac{\thet^M_X}{d^M_X\thet^M_X}
	=p^*\biggl(\frac{\thet^B_X}{d^B_X\thet^B_X}\biggr)
		+p^*\biggl(\frac{\delta_{B_X}}{e_X(\mathcal N_{B_X,B})}\biggr)
			\,\frac{\thet^M_X}{d^M_X\thet^M_X}$$
  as $L_1$-currents for any~$X\in\frg$ modulo $d_X$-exact currents.
\end{Theorem}

The first term on the right hand side gives the horizontal part of the
Bott localisation defect.
The second term is a localisation of the localisation defect to the fibres
over the zero set~$B_X$ of~$X_B$ on~$B$.
Note that the Killing vector field~$X_M$ on~$M$
is tangential to these fibres.
We can rephrase Theorem~\ref{FibreCurrentTheorem}\ as follows.
For all compactly supported continuous forms~$\alpha\in\Omega^*(M)$,
we have
	$$\int_M\frac{\thet^M_X}{d^M_X\thet^M_X}\,\alpha
	=\int_B\frac{\thet^B_X}{d^B_X\thet^B_X}\,\int_{M/B}\alpha
		+\int_B\frac{\delta_{B_X}}{e_X(\mathcal N_{B_X,B})}
			\int_{p^{-1}(B_X)/B_X}
			\frac{\thet^M_X}{d^M_X\thet^M_X}\,\alpha\;.$$
If~$X_B$ is not identically zero,
this reduces the computation of the left hand side
to the computation of integrals over manifolds of smaller dimension.

\begin{proof}
By Proposition~\ref{LeinsProp},
	$$\frac{\thet^M_X}{d^M_X\thet^M_X}\,\alpha_X
	=-\int_0^\infty\thet^M_X
	   \,e^{\textstyle t\,d^M_X\thet^M_X}
	   \,\alpha_X\,dt\;.$$
Using the vertical tangent bundle~$TF$
and the equivariant horizontal distribution~$\mathcal H=(TF)^\perp\subset TM$,
we decompose~$\thet_X=\thet_X^h+\thet_X^v
\in\Gamma(\mathcal H^*)\oplus\Gamma(T^*F)$.
If~$G$ maps fibres to fibres,
then clearly
	$$\thet_X^h
	=\frac1{2\pi i}\,\<p_*X_M,p_*\punkt\>
	=\frac1{2\pi i}\,p^*\<X_B,\punkt\>
	=p^*\thet_X^B\;,$$
hence~$\thet^M_X=p^*\thet_X^B+\thet_X^v$.

Regard the manifold~$\overline M=M\times(0,\infty)^2$,
where the action of the Killing field~$X$ is extended trivially
to~$\overline M$.
Then the form
\begin{equation*}
  \beta_X
  =e^{\textstyle d^{\overline M}_X(s\,p^*\thet_X^B+t\,\thet_X^v)}
  =\bigl(1-p^*\thet_X^B\,ds\bigr)\bigl(1-\thet_X^v\,dt\bigr)
	\,e^{\textstyle s\,p^*d^B_X\thet_X^B+t\,d^M_X\thet_X^v}
\end{equation*}
is closed on~$\overline M$.
For~$0<T\le S$ consider the domain
	$$\Omega
	=\bigl\{\,(s,t)\bigm|0\le t\le T\text{ and }t\le s\le S\,\bigr\}\;,$$
and let~$\Gamma=\Gamma_1\cup\Gamma_2\cup\Gamma_3\cup\Gamma_4$
denote the contour~$\del\Omega$,
where~$s=t$ on~$\Gamma_1$, $t=0$ on~$\Gamma_2$, $s=S$ on~$\Gamma_3$,
and~$t=T$ on~$\Gamma_4$.
It follows from the equivariant Stokes theorem that integrating~$\beta_X$
over~$\Gamma$ produces a $d_X$-exact current on~$M$.
We will consider the limits~$S\to\infty$ and~$T\to\infty$ in that order
to prove Theorem~\ref{FibreCurrentTheorem}.


Let~$I_j^0=\int_{\Gamma_j}\beta_X$ for~$j=1$, \dots, $4$,
then clearly
\begin{align}
  \begin{split}\label{IeinsFormel}
    I_1^0
    &=\int_{\Gamma_1}\beta_X
    =\int_0^T\thet^M_X\,e^{\textstyle t\,d^M_X\thet^M_X}\,dt\;,\\
    I_1^1
    &=\lim_{S\to\infty}I_1^0=I_1^0\;,\\
	\text{and}\qquad
    I_1^2
    &=\lim_{T\to\infty}I_1^1
    =-\frac{\thet_X^M}{d_X^M\thet_X^M}
  \end{split}
\end{align}
by Proposition~\ref{LeinsProp}.
Similarly,
\begin{align}
  \begin{split}\label{IzweiFormel}
    I_2^0
    &=\int_{\Gamma_2}\beta_X
    =-\int_0^Sp^*\thet^B_X\,e^{\textstyle s\,p^*d^B_X\thet^B_X}\,ds\;,\\
    I_2^1
    &=\lim_{S\to\infty}I_2^0
    =p^*\biggl(\frac{\thet_X^B}{d_X^B\thet_X^B}\biggr)\;,\\
	\text{and}\qquad
    I_2^2
    &=\lim_{T\to\infty}I_2^1=I_2^1\;.
  \end{split}
\end{align}

For the term~$I_3^0$,
we use that~$e^{S\,d^B_X\thet^B_X}$ forces localisation
on the zero set~$B_X$ of~$X_B$ in the limit~$S\to\infty$, see~\cite{BGV}.
This gives
\begin{align}
  \begin{split}\label{IdreiFormel}
    I_3^0
    &=\int_{\Gamma_3}\beta_X
    =-\int_0^T\thet^v_X
	\,e^{\textstyle S\,p^*d^B_X\thet^B_X+t\,d_M\thet_X^v}\,dt\;,\\
    I_3^1
    &=\lim_{S\to\infty}I_3^0
    =-p^*\biggl(\frac{\delta_{B_X}}{e_X\left(\mathcal N_{B_X,B}\right)}\biggr)
	\,\int_0^T\thet^M_X\,e^{\textstyle t\,d_M\thet_X^M}\,dt\;,\\
	\text{and}\qquad
    I_3^2
    &=\lim_{T\to\infty}I_3^1
    =p^*\biggl(\frac{\delta_{B_X}}{e_X\left(\mathcal N_{B_X,B}\right)}\biggr)
	\,\frac{\thet_X^M}{d_X^M\thet_X^M}\;.
  \end{split}
\end{align}

The remaining term and its limit as~$S\to\infty$ are easily computed as
\begin{align}
  \begin{split}\label{IvierFormel}
    I_4^0
    &=\int_{\Gamma_4}\beta_X
    =e^{\textstyle T\,d_X^M\thet_X^M}
	\int_0^{S-T}p^*\thet_X^B\,e^{\textstyle s\,p^*d^B_X\thet^B_X}\,ds\;,\\
	\text{and}\qquad
    I_4^1
    &=\lim_{S\to\infty}I_4^0
    =-p^*\biggl(\frac{\thet_X^B}{d_X^B\thet_X^B}\biggr)
	\,e^{\textstyle T\,d_X^M\thet_X^M}\;.
  \end{split}
\end{align}
Then as in~\eqref{LeinsZweiEst} above,
the current~$I_4^1$ localises near~$M_X$ as~$T\to\infty$,
and the limit~$I_4^2$ of~$I_4^1$ as~$T\to\infty$
can be computed using a rescaling argument.

Let~$\mathcal N=\mathcal N_{M_X,M}$ denote the normal bundle to~$M_X$ in~$M$.
We still identify the tubular neighbourhood~$U_t$ of~$M_X$
of radius~$t^{-\frac1{4}}$ for~$t=T\gg 0$
with the neighbourhood~$\mathcal N_t$ of radius~$t^{\frac14}$
by the rescaled normal exponential map~$\phy_t\colon\mathcal N_t\to U_t$
of~\eqref{PhyTDef}.
Define~$A\in\Gamma(\End\mathcal N)$ as above,
then~\eqref{ThetaXLimit} holds unchanged.
Similarly,
let~$A_0\in\Gamma(\End\mathcal N)$ denote the horizontal lift of the action
of~$\mu_X^{\mathcal N_{B_X/B}}$ on~$\mathcal N_{B_X/B}$,
then the analogue of~\eqref{ThetaXLimit} implies that
\begin{gather*}
	t\,\phy_t^*\,p^*\,X_B^*\to\<A_0\mathcal R,\punkt\>\;,
		\qquad
	t\,\phy_t^*\,p^*\,dX_B^*\to \<2A_0\punkt,\punkt\>\;,\\
		\text{and}\qquad
	t\,\phy_t^*\,p^*\norm{X_B}^2\to\norm{A_0\mathcal R}^2
\end{gather*}
uniformly on~$\mathcal N_t$.

In analogy with~\eqref{LeinsLimit},
we calculate for any continuous form~$\alpha$ that
\begin{multline*}
  \lim_{t\to\infty}\int_M
	p^*\biggl(\frac{\thet_X^B}{d_X^B\thet_X^B}\biggr)
	\,e^{\textstyle t\,d_X^M\thet_X^M}\,\alpha\\
  \begin{aligned}
    &=\lim_{t\to\infty}\int_{\mathcal N_t}\sum_{j,k}
	\,\phy_t^*p^*\biggl(\frac{t\,X_B^*\,(t\,dX_B^*)^j}
		{(2\pi i\,t\norm{X_B}^2)^{j+1}}\biggr)
	\,\frac{(t\,\phy_t^*dX_M^*)^k}{(2\pi i)^kk!}
	\,e^{-t\,\phy_t^*\norm{X_M}^2}
	\,\phy_t^*\alpha\\
    &=\lim_{t\to\infty}
      \int_{\mathcal N}\sum_{j,k}
	\,\frac{\<A_0\mathcal R,\punkt\>\,\<2A_0\punkt,\punkt\>^j
		\,\<2A\punkt,\punkt\>^k}
		{(2\pi i)^{j+k+1}\,\norm{A_0\mathcal R}^{2j+2}\,k!}
	\,e^{-\norm{A\mathcal R}^2}\,\pi^*\alpha\\
    &=0\;.\\
  \end{aligned}
\end{multline*}
With~\eqref{IvierFormel}, this implies that
	$$I_4^2=\lim_{T\to\infty}I_4^1=0\;.$$
Theorem~\ref{FibreCurrentTheorem}\
follows from~\eqref{IeinsFormel}--\eqref{IdreiFormel},
because~$X_1^2+X_2^2+X_3^2+X_4^2$ is exact.
\end{proof}

\subsection{The Bott localisation defect on \texorpdfstring{$G/H$}{G/H}}\label{HomogenKapitel}
Let~$G\supset K\supset S$ be compact Lie groups,
then
\begin{equation}\label{HomFibDiagramm}
  \begin{CD}
	P=K/S@>>>E=G/S\\
	&&@VVpV\\
	&&Q=G/K
  \end{CD}
\end{equation}
is a $G$-equivariant fibration.
The left and the right hand triangle in the following diagram
represent two such equivariant fibrations with~$K=T$ and~$K=H$, respectively.
\begin{equation}\label{GruppenDiagramm}
  \begin{matrix}
    &&\!\!G\!\!\\
    &\diagup&&\diagdown\\
    T\!\!\!\!&&\!\!\smash{\Bigg|}\!\!&&\!\!\!\!H\\
    &\diagdown&&\diagup\\
    &&\!\!S\!\!
  \end{matrix}
\end{equation}
We apply Theorem~\ref{FibreCurrentTheorem}\ to these fibrations
to compute
	$$\int_{G/H}\frac{\thet_X}{d_X\thet_X}
		\,\Adach_X\bigl(T(G/H),\nabla\LC\bigr)
		\,\ch_X\bigl(\mathcal E/\mathcal S,\nabla^{\mathcal E}\bigr)\;.$$

Together with the formula for~$\eta_{e^{-X}}(\tilde D)+h_{e^{-X}}(\tilde D)$
in~\cite{G1}, \cite{G2},
we obtain a formula for~$\eta(D)+h(D)$ up to a possible contribution in~$2\Z$
coming from the spectral flow.

Recall the definition of~$X_E$ for~$X\in\frg$ in~\eqref{XMDef}.
With the notation of~\eqref{TMiso},
we represent~$X_E$ by~$\hat X\colon G\to\fre=\frs^\perp$ with
	$$\hat X_E(g)=-\bigl(\Ad_g^{-1}X\bigr)_\fre\;.$$
The Lie derivative~$\mathcal L_X^\kappa$ acts on~$V^\kappa E$ by
	$$\widehat{\mathcal L_X^\kappa s}
	=\frac d{dt}\Bigr|_{t=0}\hat s\bigl(e^{-tX}g\bigr)
	=-\bigl(\Ad_g^{-1}X\bigr)\bigl(\hat s\bigr)\;,$$
so the moment of~$X$ with respect to the reductive connection is simply
\begin{equation}\label{MomentNullFormel}
  \widehat{\mu^{0,\kappa}_Xs}
  =\widehat{\mathcal L_X^\kappa s}-\widehat{\nabla^{0,\kappa}_{X_E}s}
  =-\bigl(\Ad_g^{-1}X\bigr)_\frh\bigl(\hat s\bigr)
  =\kappa_{*\left(\Ad_g^{-1}X\right)_\frh}\hat s\;.
\end{equation}
The equivariant curvature of~$\nabla^{0,\kappa}$ is then given by
\begin{equation}\label{FNullXFormel}
  F^{0,\kappa}_X
  =\bigl(\nabla^{0,\kappa}-2\pi i\,\iota_{X_E}\bigr)^2
	+2\pi i\,\mathcal L_X^\kappa
  =F^{0,\kappa}+2\pi i\,\mu^{0,\kappa}_X
  =\kappa_{*\left(-[\punkt,\punkt]_\frh+2\pi i\,(\Ad_g^{-1}X)_\frh\right)}\;.
\end{equation}

We will assume that~$\mathcal E=\mathcal S^\kappa E\to E$ is constructed
as in~\eqref{SkappaDef}.
If we regard the Dirac operator~$D^\kappa=D^{\frac1{2},\kappa}$
on~$\mathcal S^\kappa E$,
then the equivariant twisting curvature in~\eqref{FESDef}
is just the equivariant curvature of the reductive connection
on~$V^\kappa E$
(even if this bundle does not exist globally),
so formally
\begin{multline}\label{ChChFormel}
  \ch_X\bigl(\mathcal E/\mathcal S,\nabla^{\mathcal E}\bigr)
  =\ch_X\bigl(V^\kappa E,\nabla^0\bigr)\\
  =\trace_{V^\kappa}\biggl(e^{\smash{\textstyle\frac{\kappa_{*[\punkt,\punkt]_\frh}}{2\pi i}}
	-\kappa_{\smash{*(\Ad_g^{-1}X)_\frh}}}\biggr)
  =\chi_H^\kappa
    \biggl(e^{{\textstyle\frac{[\punkt,\punkt]_\frh}{2\pi i}}
	-(\Ad_g^{-1}X)_\frh}\biggr)\;.
\end{multline}

Let~$P\hookrightarrow E\to Q$ be an equivariant fibration
as in~\eqref{HomFibDiagramm} with~$K=H$,
so the fibre~$P=H/S$ is a flag manifold.
An $\Ad_G$-invariant metric on~$\frg$ induces a normal metric on~$E$
such that~$p$ is an equivariant Riemannian submersion.
We write~$\frp=\frs^\perp\cap\frh$ and~$\frq=\frh^\perp$,
so~$\fre=\frp\oplus\frq$.
By $\Ad_G$-invariance,
the isotropy representation of~$E$ splits as~$\pi=\phy\oplus\iota^*\psi$,
where~$\phy$ and~$\psi$ are the isotropy representations of~$P$ and~$Q$,
and~$\iota\colon S\to H$ is the inclusion.
Thus the tangent bundle~$T E$ of the total space
splits naturally as~$p^*TQ\oplus TP$
with vertical tangent bundle~$TP=G\times_S\frp$
and horizontal complement~$p^*TQ=G\times_S\frq$.

We define~$\nabla^{\lambda,P}$, $\nabla^{\mu,Q}$ as in~\eqref{NablaTDef} above,
regarding~$\nabla^{\lambda,P}$ as a connection
on the vertical tangent bundle~$TP\to E$.
Then
	$$\widehat{\nabla^{\lambda,P}_UY}
	=\hat U(\hat Y)+\lambda\,[\hat U,\hat Y]_\frp
		\qquad\text{and}\qquad
	\widehat{\nabla^{\mu,Q}_VW}
	=\hat V(\hat W)+\mu\,[\hat V,\hat W]_\frq$$
for~$S$-equivariant~$\hat U\colon G\to\fre$, $\hat Y\colon G\to\frp$
and~$H$-equivariant~$\hat V$, $\hat W\colon G\to\frq$.

Let~$\kappa\in\frs$ be a weight of~$\frh$,
and consider the bundle~$V^{\kappa+\rho_H}P\to P$
associated to the $S$-representation with highest weight~$\kappa+\rho_H$.
Note that this bundle lives on the universal cover of~$P=H/S$.
However,
its equivariant Chern form~$\ch_X(V^{\kappa+\rho_H}P,\nabla^0\bigr)$
descends to~$P$,
even if~$P$ is not simply connected.
Similarly, the character~$\chi^\kappa_H(e^X)$
is well-defined in terms of~$X\in\frh$ even if~$\chi^\kappa_H$
is maybe not defined as a function on~$H$.
Bott's localisation formula and the equivariant index theorem imply that
\begin{equation}\label{KappaIndThm}
  \chi^\kappa_H\bigl(e^{-X}\bigr)
  =\int_{P}\Adach_X\bigl(TP,\nabla^{\lambda,P}\bigr)
	\,\ch_X\bigl(V^{\kappa+\rho_H}P,\nabla^0\bigr)\;,
\end{equation}
see~\cite{BGV}, Section~8.2.

The curvature~$F^{\lambda,P}$ and the moment~$\mu_X^{\lambda,P}$ on~$E$
are given by
\begin{align*}
  \widehat{F^{\lambda,P}_{V,W}Y}
  &=-\phy_{*[\hat V,\hat W]_\frs}\hat Y
	-\lambda\,\bigl[[\hat V,\hat W]_\frp,\hat Y\bigr]_\frp
	+\lambda^2\,\bigl[\hat V,[\hat W,\hat Y]_\frp\bigr]_\frp
	-\lambda^2\,\bigl[\hat W,[\hat V,\hat Y]_\frp\bigr]_\frp\\
  \text{and}\qquad
  \widehat{\mu_X^{\lambda,P}Y}
  &=\phy_{*(\Ad_h^{-1}X)_\frs}\hat 
	Y+\lambda\,\bigl[(\Ad_h^{-1}X)_\frp,\hat y\bigr]_\frp\;.
\end{align*}
Using these equations,
one checks that
the equivariant curvature of the connection~$\nabla^{\lambda,P}$
and the Killing field~$X$ on~$E$
is formally the same as the equivariant curvature of~$\nabla^{\lambda,P}$
and the Killing
field~$-\Ad_h{\textstyle\frac{[\punkt_\frq,\punkt_\frq]_\frh}{2\pi i}}+X$,
considered only on~$P$.
By~\eqref{ChChFormel} and~\eqref{KappaIndThm},
this gives
\begin{multline}\label{CWKirillovFormel}
  \ch_X\bigl(V^{\kappa}Q,\nabla^0\bigr)
    =\int_{P}\Adach_{\smash
	{-\Ad_h{\textstyle\frac{[\punkt_\frq,\punkt_\frq]_\frh}{2\pi i}}+X}}
		\bigl(TP,\nabla^{\lambda,P}\bigr)
      \,\ch_{\smash
	{-\Ad_h{\textstyle\frac{[\punkt_\frq,\punkt_\frq]_\frh}{2\pi i}}+X}}
		\bigl(V^{\kappa+\rho_H}P,\nabla^0\bigr)\\
    =\int_{ E/Q}\Adach_X\bigl(TP,\nabla^{\lambda,P}\bigr)
	\,\ch_X\bigl(V^{\kappa+\rho_H} E,\nabla^0\bigr)\;,
\end{multline}
where we integrate over the fibres of~$E\to Q$ in the last line.

Let~$p^*\nabla^{\mu,Q}\oplus\nabla^{\lambda,P}$ denote the product
connection on~$T E=p^*TQ\oplus TF$.
By multiplicativity of the $\Adach$-form,
we have
	$$\Adach_X\bigl(T E,p^*\nabla^{\mu,Q}\oplus\nabla^{\lambda,P}\bigr)
	=p^*\Adach_X\bigl(TQ,\nabla^{\mu,Q}\bigr)
		\,\Adach_X\bigl(TP,\nabla^{\lambda,P}\bigr)\;.$$
Now, Theorem~\ref{FibreCurrentTheorem}\ and equation~\eqref{CWKirillovFormel}
have the following implication.

\begin{Proposition}\label{RightProposition}
  Let~$G\supset H$ be compact Lie groups,
  and let~$S\subset H$ be a maximal torus.
  If the Killing field~$X_{Q}$ has no zeros,
  then for any~$\lambda$ and~$\mu$,
  \begin{multline*}
	\frac{\thet_X}{d_X\thet_X}\,\Adach_X\bigl(TQ,\nabla^{\mu,Q}\bigr)
		\,\ch_X\bigl(V^\kappa Q,\nabla^0\bigr)\\
	=\int_{ E/Q}\frac{\thet_X}{d_X\thet_X}
		\,\Adach_X\bigl(T E,p^*\nabla^{\mu,Q}\oplus\nabla^{\lambda,P}\bigr)
		\,\ch_X\bigl(V^{\kappa+\rho_H} E,\nabla^0\bigr)\;.\quad\qed
  \end{multline*}
\end{Proposition}

To evaluate the right hand side of Proposition~\ref{RightProposition},
we now consider the left triangle in~\eqref{GruppenDiagramm}.
In particular, $T\subset G$ is a maximal torus containing~$S$.
We assume in addition that~$X$ is regular.
Then the only fixed points of~$X$ on~$G/T$ are the isolated fixpoints~$wT/T$.
In particular,
the normal bundle to~$N_G(T)/T$ in~$G/T$ is just the tangent bundle~$T(G/T)$,
and its equivariant Euler form is invertible near~$N_G(T)/T$.
We find by Proposition~\ref{RightProposition} that
\begin{multline*}
  \int_{Q}\frac{\thet_X}{d_X\thet_X}\,\Adach_X\bigl(TQ,\nabla^{\mu,Q}\bigr)
	\,\ch_X\bigl(V^\kappa Q,\nabla^0\bigr)\\
  =\int_{E}\Adachsl_X\bigl(T E,\nabla^0,
	p^*\nabla^{\mu,Q}\oplus\nabla^{\lambda,P}\bigr)
	\,\ch_X\bigl(V^{\kappa+\rho_H} E,\nabla^0\bigr)\\
  +\int_{E}\frac{\thet_X}{d_X\thet_X}
	\,\Adach_X\bigl(T E,\nabla^0\bigr)
	\,\ch_X\bigl(V^{\kappa+\rho_H} E,\nabla^0\bigr)\;.
\end{multline*}
Applying Theorem~\ref{FibreCurrentTheorem} once more gives
\begin{multline}\label{TotalFormelEins}
  \int_{Q}\frac{\thet_X}{d_X\thet_X}\,\Adach_X\bigl(TQ,\nabla^{\mu,Q}\bigr)
	\,\ch_X\bigl(V^\kappa Q,\nabla^0\bigr)\\
  \begin{aligned}
	&=\int_{E}\Adachsl_X\bigl(T E,\nabla^0,
		p^*\nabla^{\mu,Q}\oplus\nabla^{\lambda,P}\bigr)
		\,\ch_X\bigl(V^{\kappa+\rho_H} E,\nabla^0\bigr)\\
	&\qquad
	+\int_{G/T}\frac{\thet_X}{d_X\thet_X}
		\int_{ E/(G/T)}\Adach_X\bigl(T E,\nabla^0\bigr)
			\,\ch_X\bigl(V^{\kappa+\rho_H} E,\nabla^0\bigr)\\
	&\qquad+\int_{G/T}\frac{\delta_{N_G(T)/T}}{e_X\bigl(T(G/T)\bigr)}
		\int_{N_G(T)/S}\frac{\thet_X}{d_X\thet_X}
		\,\Adach_X\bigl(T E,\nabla^0\bigr)
			\,\ch_X\bigl(V^{\kappa+\rho_H} E,\nabla^0\bigr)
  \end{aligned}
\end{multline}

We can still simplify this expression.
First of all,
the reductive connection respects the splitting~$TE\cong p^*TQ\oplus TP$.
We can thus fix~$\lambda=0$
and rewrite the first term of the right hand side of~\eqref{TotalFormelEins}\
as
\begin{multline*}
  \int_{E}\Adachsl_X\bigl(T E,\nabla^0,
	p^*\nabla^{\mu,Q}\oplus\nabla^{0,P}\bigr)
	\,\ch_X\bigl(V^{\kappa+\rho_H} E,\nabla^0\bigr)\\
  =\int_{E}\Adachsl_X\bigl(p^*TQ,\nabla^0,p^*\nabla^{\mu,Q}\bigr)
	\,\Adach\bigl(TP,\nabla^{0,P}\bigr)
	\,\ch_X\bigl(V^{\kappa+\rho_H} E,\nabla^0\bigr)\;.
\end{multline*}
We want to show that for~$\mu=0$,
the expression above vanishes.

Note that the curvature of the reductive connection on an equivariant vector
bundle over~$E$
depends only on the $(G\times_H\frs)$-valued two form~$[\punkt,\punkt]_\frs$
by~\eqref{CurvatureNullFormel}.
However, this form vanishes on~$\frg/\frh\otimes\frh/\frs$,
so we find that
\begin{equation}\label{IntegrandParityFormel}
  \Adach_X\bigl(TP,\nabla^{0,P}\bigr)
	\,\ch_X\bigl(V^{\kappa+\rho_H} E,\nabla^0\bigr)
  \in\Gamma\bigl(p^*\Lambda^\even T^*Q
	\otimes\Lambda^\even T^*P\bigr)\;.
\end{equation}

It remains to analyse the Chern-Simons
class~$\Adachsl_X(p^*TQ,\nabla^0,p^*\nabla^{0,Q})$.
Once again,
let~$\frp=\frs^\perp\cap\frh$ and~$\frq=\frh^\perp\subset\frg$.
Note that a vector field~$V$ on~$Q$
is given by an $H$-equivariant function~$\hat V\colon G\to\frq$.
This function also describes the horizontal lift of~$V$ to~$p^*TQ\subset TE$.
Therefore,
the pull-back connection~$p^*\nabla^{0,Q}$ on
the horizontal tangent bundle~$p^*TQ\subset TE$
is given by
\begin{equation}\label{PsternNablaFormel}
  \widehat{p^*\nabla^{\mu,Q}_VW}
  =\hat V\bigl(\hat W\bigr)
	+\bigl[\hat V_\frp,\hat W\bigr]
\end{equation}
for~$S$-equivariant functions~$V\colon G\to\fre$ and~$W\colon G\to\frq$,
since then
	$$\widehat{p^*\nabla^{\mu,Q}_VW}
	=\hat V_\frq\bigl(\hat W\bigr)$$
for an $H$-equivariant~$\hat W$.

We fix a family of connections~$\nabla^{0,\lambda}$
on~$p^*TQ\subset T E$ by
	$$\widehat{\nabla^{0,\lambda}_VW}
	=\hat V(\hat W)+\lambda\bigl[\hat V_\frp,\hat W\bigr]\;.$$
By~\eqref{NablaNullDef}\ and~\eqref{PsternNablaFormel},
we have~$\nabla^{0,0}=\nabla^0$ and~$\nabla^{0,1}=p^*\nabla^{0,Q}$.
A straightforward computation gives the differential and the curvature
of this family of connections as
\begin{alignat*}2
	\frac\del{\del\lambda}\widehat{\nabla^{0,\lambda}_V}
	&=\bigl[\hat V_\frp,\punkt\bigr]
	&\quad&\in\Lambda^1\frp\otimes\End\frq\;,\\
\noalign{\smallskip}
		\text{and}\qquad
	\widehat{F^{0,\lambda}_{V,W}}
	&=-\psi_{*\left[\hat V,\hat W\right]_\frs}
		-\lambda\,\bigl[[\hat V,\hat W]_\frp,\punkt\bigr]\\
	&\qquad
		+\lambda^2\,\bigl[\hat V_\frp
			,[\hat W_\frp,\punkt]\bigr]
		-\lambda^2\,\bigl[\hat W_\frp
			,[\hat V_\frp,\punkt]\bigr]
	&\quad&\in\bigl(\Lambda^2\frp\oplus\Lambda^2\frq\bigr)\otimes\End\frq\;,
\end{alignat*}
where we have used that~$[\frp,\frq]\subset\frq$
and~$\frq\cap\frs=0$.
This implies that
\begin{equation}\label{ChernSimonsParityFormel}
  \Adachsl_X\bigl(p^*TQ,\nabla^0,p^*\nabla^{0,Q}\bigr)
  \in\Gamma\bigl(p^*\Lambda^\even T^*Q\otimes\Lambda^\odd T^*P\bigr)\;.
\end{equation}
Combining~\eqref{IntegrandParityFormel} and~\eqref{ChernSimonsParityFormel}
with the fact that~$\dim Q$ is odd while~$\dim P$ is even,
we see that
\begin{multline*}
  \int_{E}\Adachsl_X\bigl(T E,\nabla^0,
	p^*\nabla^0\oplus\nabla^0\bigr)
	\,\ch_X\bigl(V^{\kappa+\rho_H} E,\nabla^0\bigr)\\
  =\int_{E}\Adachsl_X\bigl(p^*TQ,\nabla^0,p^*\nabla^0\bigr)
	\,\Adach_X\bigl(TP,\nabla^0\bigr)
	\,\ch_X\bigl(V^{\kappa+\rho_H} E,\nabla^0\bigr)
  =0\;.
\end{multline*}
Of course,
this conclusion would fail for most other possible choices of connections.

The analogue of~\eqref{IntegrandParityFormel} for the left hand side
of \eqref{GruppenDiagramm} is
	$$\Adach_X\bigl(T E,\nabla^0\bigr)
		\,\ch_X\bigl(V^{\kappa+\rho_H} E,\nabla^0\bigr)
	\in\Gamma\bigl(p^*\Lambda^\even T^*(G/T)
		\otimes\Lambda^\even T^*(T/S)\bigr)\;.$$
It implies that
the middle term of the right hand side of~\eqref{TotalFormelEins} vanishes.
Note that this conclusion would also fail for many other possible choices
of connections.

Finally,
over~$W_G=N_G(T)/T\subset G/T$,
the Killing field~$X_{E}$ becomes tangential to the fibres.
Let~$w=nT\in N_G(T)/T=W_G$,
then~$X_{E}|_{nT}$ is diffeomorphic to the Killing field~$wX_{T/S}$
on~$T/S$.
Let us summarise our computations so far.

\begin{Proposition}\label{TotalProposition}
  Let~$G\supset H\supset S$ and~$G\supset T\supset S$ be as above,
  and assume that~$X\in\frt$ is regular
  and that the Killing field~$X_{G/H}$ has no zeros.
  Then
  \begin{multline*}
    \int_{G/H}\frac{\thet_X}{d_X\thet_X}\,\Adach_X\bigl(T(G/H),\nabla^0\bigr)
	\,\ch_X\bigl(V^\kappa(G/H),\nabla^0\bigr)\\
    =\sum_{w\in W_G}\frac{1}{e_{wX}\bigl(T(G/T)\bigr)}
      \int_{T/S}\frac{\thet_{wX}}{d_{wX}\thet_{wX}}
	\,\Adach_{wX}\bigl(T(G/S),\nabla^0\bigr)
	\,\ch_{wX}\bigl(V^{\kappa+\rho_H}(G/S),\nabla^0\bigr)\;.
  \end{multline*}
\end{Proposition}

\subsection{Evaluation of the Bott localisation defect}\label{EvalAbschnitt}
In the previous subsection,
we have reduced the Bott localisation defect
to the quotient of the maximal tori.
We will now give two formulas in representation theoretic terms.

We start with a Lemma concerning the Weyl groups of~$G$ and~$H$.

\begin{Lemma}\label{SoergelLemma}
  Let~$H\subset G$ be a pair of compact Lie groups
  with maximal tori~$S\subset T$,
  and let~$W_G$, $W_H$ be the corresponding Weyl groups.
  Then~$W_H$ is a subgroup of
	$$\bigl\{\,w|_{\frs}\bigm|w\in W_G\text{ and }
		w(\frs)=\frs\,\bigr\}\subset\Aut(S)\;.$$
\end{Lemma}

\begin{proof}
Let~$w'=n'S\in W_H=N_H(S)/S$,
then~$w'$ acts on~$S$ by~$s\mapsto n'sn^{\prime-1}$.
Clearly, $n'Tn^{\prime-1}\subset Z_G(S)$ is a maximal torus
in the centraliser of~$S$ in~$G$.
In particular, we find~$z\in Z_G(S)$ such that~$zn'T(zn')^{-1}=T$,
and~$zn's(zn')^{-1}=w'(s)$ for all~$s\in S$.
Thus, $w=zn'T\in W_G=N_G(T)/T$ acts on the subset~$S\subset T$ as~$w'$.
\end{proof}

To state our main result,
let~$\delta\in i\frt^*$ be defined as in~\eqref{DeltaDef}.
Let~$\prod_{\beta\in W_\pi^+}$ denote the product over all positive weights
of the isotropy representation~$\pi$,
counted with multiplicity by abuse of notation.
Let~$\Adach\colon z\to\frac{z/2}{\sinh(z/2)}$ denote the $\Adach$-function.
We will denote by~$X|_\frs$ the orthogonal projection of~$X\in\frt$
onto~$\frs$.
If~$f\colon\frt\to\C$ satisfies~$f(X)=f(wX)$
for all~$w\in W_G$ that map~$\frs$ to itself,
then we set
\begin{equation}\label{SGHDef}
  S_{G/H}(f)(X)
  =\frac1{\#W_H}\,\sum_{w\in W_G}\sign(w)\,f(wX)\;.
\end{equation}
If~$W_H$ can be identified with a subgroup of~$W_G$,
this amounts to summing over~$W_G/W_H$.

\begin{Theorem}\label{BLDThm}
  Let~$M=G/H$ be a quotient of compact Lie groups.
  If~$\rk G=\rk H+1$,
  then
  \begin{align*}
	\int_M\frac{\thet_X}{d_X\thet_X}\,\Adach_X&\bigl(TM,\nabla^0\bigr)
		\,\ch_X\bigl(V^\kappa M,\nabla^0\bigr)\\
	&=A_G\Biggl(\frac{e^{(\kappa+\rho_H)(\punkt|_\frs)}}{\delta(\punkt)}
			\,\prod_{\beta\in\Delta_G^+}
			\Adach\bigl(\beta(\punkt|_\frs)\bigr)\Biggr)(-X)
		\,\prod_{\beta\in\Delta_\frg^+}\frac{-1}{\beta(X)}
		\tag1\\
	&=S_{G/H}\Biggl(\frac{\chi_H^\kappa\bigl(e^{\punkt|_\frs}\bigr)}
				{\delta(\punkt)}
			\,\prod_{\beta\in\Delta_G^+}
				\frac{\beta(\punkt|_\frs)}{\beta(\punkt)}
			\,\prod_{\beta\in W_\pi^+}
				\frac1{2\sinh\bigl(\frac\beta2(\punkt|_\frs)
					\bigr)}\Biggr)(-X)\;.
		\tag2
  \end{align*}
  Otherwise,
  \begin{align*}
    \int_M\frac{\thet_X}{d_X\thet_X}\,\Adach_X\bigl(TM,\nabla^0\bigr)
		\,\ch_X\bigl(V^\kappa M,\nabla^0\bigr)
    &=0\;.\tag3
  \end{align*}
\end{Theorem}

Formula~(2) gives an advantage for explicit computations
if the subgroup~$H$ has a large Weyl group
and the character~$\chi_H^\kappa$ is known,
e.g., if~$D^\kappa=D$ is the untwisted Dirac operator.
We have used this formula in Section~\ref{SymmAbschnitt}
when dealing with spheres and odd Grassmannians.
In some cases,
one can even improve on~(2).
Suppose that there exists another subgroup~$K\subset G$,
such that~$H\subset K$ and~$S$ is a maximal torus of~$K$.
Then one obtains a similar formula where~$K$ replaces~$H$,
$\pi$ becomes the isotropy representation of~$G/K$,
and~$\kappa$ gets replaced by~$\kappa+\rho_H-\rho_K$.

\begin{proof}
  We evaluate the right hand side of Proposition~\ref{TotalProposition}
  term by term.
  The tangent bundle~$T(G/T)$ is oriented by the choice of~$\Delta_G^+$
  as in~\eqref{GTOrientation},
  so for~$w=nT$,
	$$\frac{\delta_{nT/T}}{e_X\bigl(T(G/T)\bigr)}
	=\prod_{\beta\in\Delta_G^+}\frac{-1}{\beta(wX)}
	=\sign(w)\prod_{\beta\in\Delta_G^+}\frac{-1}{\beta(X)}\;,$$
  independent of any connection on~$T(G/T)$,
  since we evaluate at isolated fixpoints.
  
  Because the vertical tangent bundle of~$G/S\to G/T$
  is trivial, $G$-invariant and parallel with respect to~$\nabla^0$,
  the equivariant $\Adach$-form of~$T(G/S)$ is given as
	$$\Adach_{wX}\bigl(T(G/S),\nabla^0\bigr)|_{T/S}
	=\Adach_{wX}\bigl(p^*T(G/T),\nabla^0\bigr)|_{T/S}
	=\prod_{\beta\in\Delta_G^+}\Adach\bigl(\beta(-wX|_\frs)\bigr)$$
  by~\eqref{MomentNullFormel};
  because~$T/S$ is one-dimensional,
  only the moment~$\mu^0$ of~$\nabla^0$ enters.
  Similarly,
  the equivariant Chern character form equals
	$$\ch_{wX}\bigl(V^{\kappa+\rho_H}(G/S),\nabla^0\bigr)|_{T/S}
	=e^{(\kappa+\rho_H)(-wX|_\frs)}\;.$$

  Assume that~$\rk G=\rk H+1$.
  The expression~$\int_{T/S}\frac{\thet_{wX}}{d_{wX}\thet_{wX}}$
  is clearly independent of the metric chosen,
  so we may assume that~$\Vol(T/S)=2\pi$.
  Then for a positively oriented unit vector~$E$,
  we have~$\delta(E)=i$ and~$\delta(X)=i\<X,E\>$.
  Then
	$$\int_{T/S}\frac{\thet_{wX}}{d_{wX}\thet_{wX}}
	=-\int_{T/S}\frac{wX^*_{T/S}}{2\pi i\,\norm{wX_{T/S}}^2}
	=-\frac{\Vol(T/S)}{2\pi i\,\<wX,E\>}
	=\frac1{\delta(-wX)}\;.$$
  This proves~(1).

  To prove~(2),
  we use Lemma~\ref{SoergelLemma} to simplify~(1).
  If we assume that~$\rk G=\rk H+1$,
  there are two possibilities.
  Either,
  each element of~$W_H$ corresponds to precisely one element of~$W_G$;
  then these elements form a subgroup of~$W_G$ that we identify with~$W_H$.
  Or each element of~$W_H$ corresponds to precisely two elements of~$W_G$
  which differ by the reflexion on~$S$;
  then we can identify~$W_H$ with a subgroup of~$W_G^{\even}$,
  or with those elements that preserve the normal orientation of~$S$ in~$T$.
  We fix one of these identifications.
  Let~$\sign_G(w)=\det(w|_\frt)$ and~$\sign_H(w)=\det(w|_\frs)$
  denote the sign of an element~$w\in W_H\subset W_G$ as an element of the
  Weyl groups~$W_G$ and~$W_H$,
  respectively.
  These signs are related by
  \begin{equation}\label{WeylSignFormel}
    \sign_G(w)\,\delta(wX)=\sign_H(w)\,\delta(X)\;.
  \end{equation}

  As representations of~$\frs$, $\frg/\frt\oplus\R\,E$
  is isomorphic to~$\frh/\frs\oplus\pi$.
  If we have chosen~$\Delta_G^+$ and~$\Delta_H^+$ carefully,
  the restriction to~$\frs\subset\frt$
  maps the positive roots~$\Delta_G^+$ bijectively
  onto~$\Delta_H^+\mathrel{\dot\cup}W_\pi^+$,
  where the positive weights of~$\pi$ are counted with multiplicity.
  We assume that this is the case,
  even though our final formula does not depend on this choice,
  only on the compatibility of orientations in~\eqref{OrientFormel}.
  Recall the Weyl formulas,
  \begin{equation}\label{WeylFormulas}
    \prod_{\beta\in\Delta_H^+}2\sinh\biggl(\frac\beta2(X)\biggr)
    =A_H\bigl(e^{\rho_H}\bigr)(X)
	\qquad\text{and}\qquad
    \frac{A_H\bigl(e^{\kappa+\rho_H}\bigr)}{A_H\bigl(e^{\rho_H}\bigr)}(X)
    =\chi_H^\kappa\bigl(e^X\bigr)\;,
  \end{equation}
  where~$A_H$ denotes the alternating sum over the Weyl group of~$H$.
  Regarding~$W_H$ as a subgroup of~$W_G$ and using~\eqref{WeylSignFormel}
  and~\eqref{WeylFormulas},
  we find that
  \begin{multline*}
    A_H\Biggl(\frac{e^{(\kappa+\rho_H)(\punkt|_\frs)}}{\delta(\punkt)}
	\,\prod_{\beta\in\Delta_G^+}
		\Adach\bigl(\beta(\punkt|_\frs)\bigr)\Biggr)(-X)
	\,\prod_{\beta\in\Delta_G^+}\frac1{\beta(-X)}\\
    \begin{aligned}
      &=\frac1{\delta(-X)}
	\,\frac{A_H\bigl(e^{\kappa+\rho_H})}
		{A_H\bigl(e^{\rho_H})}\bigl(-X|_\frs\bigr)
	\prod_{\beta\in\Delta_G^+}
		\frac{\Adach\bigl(\beta(X|_\frs)\bigr)}{\beta(-X)}
	\prod_{\beta\in\Delta_H^+}
		2\sinh\biggl(\frac\beta2(X|_\frs)\biggr)\\
      &=\frac{\chi_H^\kappa\bigl(e^{-X|_\frs}\bigr)}{\delta(-X)}
	\prod_{\beta\in\Delta_G^+}\frac{\beta(-X|_\frs)}{\beta(-X)}
	\prod_{\beta\in W_\pi^+}
		\frac1{2\sinh\bigl(\frac\beta2(-X|_\frs)\bigr)}\;.
    \end{aligned}
  \end{multline*}
  Replacing~$X$ by~$wX$ and
  summing over~$wX$ for~$w\in W_G/W_H$,
  we obtain~(2).

  Finally,
  assume that~$\rk G\ne\rk H+1$.
  First of all, if~$\rk G-\rk H$ is even,
  then so is~$\dim M$,
  and the integral in Theorem~\ref{BLDThm}\ vanishes for parity reasons.
  Otherwise, $\rk G-\rk H=\dim T/S\ge 3$.
  Because~$T$ is abelian,
  the two-form part of the equivariant curvature~$F^{0,\kappa}_{wX}$
  vanishes on~$T/S$ by~\eqref{FNullXFormel}.
  This implies that~$\Adach_{wX}(T(G/S),\nabla^0)$
  and~$\ch_{wX}(V^{\kappa+\rho_H}(G/S),\nabla^0)$
  have no components of non-zero exterior degree.
  Also,
  the Killing vector field on~$T/S$ generated by~$wX\in\frt$
  is parallel,
  so~$\frac{\thet_{wX}}{d_{wX}\thet_{wX}}\in\Omega^1(T/S)$.
  Then
	$$\frac{\thet_{wX}}{d_{wX}\thet_{wX}}
		\,\Adach_{wX}\bigl(T(G/S),\nabla^0\bigr)
		\,\ch_{wX}\bigl(V^{\kappa+\rho_H}(G/S),\nabla^0\bigr)
	\in\Omega^1(T/S)\;,$$
  and the integral over~$T/S$ vanishes,
  which proves~(3).
\end{proof}

\subsection{Evaluation of \texorpdfstring{$\eta$}{eta}-invariants}\label{EqEtaSubsection}
We combine Theorem~\ref{BLDThm}\ with Theorem~\ref{GdreiThm}\ and Theorem~\ref{EtaThm}\
to establish our final formula
for the infinitesimally equivariant $\eta$-invariant
of the Dirac operator~$D^\kappa$ on~$M=G/H$.

Let~$\hat G$ denote the set of equivalence classes of irreducible
unitary representations of~$G$.
Because~$G$ is compact,
all these representations are finite dimensional.
Recall that by Frobenius reciprocity and the Peter-Weyl theorem,
	$$\Gamma(S^\kappa M)
	=\overline{\bigoplus_{\gamma\in\hat G}V^\gamma
		\oplus\Hom_H\bigl(V^\gamma,S\otimes V^\kappa\bigr)}\;.$$
Because~$D^\kappa$ and~$\Dsl^\kappa$ are $G$-equivariant,
they respect this decomposition,
and we have
	$$D^\kappa
	=\overline{\bigoplus_{\gamma\in\hat G}
		\id_{V^\gamma}\otimes{}^{\gamma\!} D^\kappa}\;,$$
and similarly for~$\Dsl^\kappa$.
Let~$\eta({}^{\gamma\!} D^\kappa)\in\Z$ and~$h({}^{\gamma\!} D^\kappa)\in\Z$
denote the $\eta$-invariant
and the dimension of the kernel of~${}^{\gamma\!} D^\kappa$
acting on~$\Hom_H(V^\gamma,S\otimes V^\kappa)$.
We also recall the definitions of~$\delta$ and~$\alpha\in\frt^*$
in~\eqref{DeltaDef}, \eqref{AlphaDef},
and the definition of~$W_\pi^+$  and~$S_{G/H}$ in section~\ref{EvalAbschnitt}.

\begin{Theorem}\label{MainTheorem}
  Assume that~$\rk G=\rk H+1$,
  then
  \begin{align*}
	\eta_X\bigl(D^\kappa\bigr)
	&=2\sum_{w\in W_G}\frac{\sign(w)}{\delta(wX)}
	\,\Biggl(\prod_{\beta\in\Delta_G^+}\Adach\bigl(\beta(wX)\bigr)
		\,\Adach\bigl(\delta(wX)\bigr)
		\,e^{-\left(\alpha-\frac\delta2\right)(wX)}\\
	&\kern10em
	-\prod_{\beta\in\Delta_G^+}\Adach\bigl(\beta(wX|_\frs)\bigr)
		\,e^{-(\kappa+\rho_H)(wX|_\frs)}\Biggr)
		\,\prod_{\beta\in\Delta_G^+}\frac{-1}{\beta(X)}\tag1\\
	&\qquad
		+2\int_M\Adachsl_X\bigl(TM,\nabla^0,\nabla\LC\bigr)
			\,\ch_X\bigl(V^\kappa M,\nabla^\kappa\bigr)\\
	&\qquad
		+\sum_{\gamma\in\hat G}\chi_G^\gamma\bigl(e^{-X}\bigr)
			\Bigl(\eta\bigl({}^{\gamma\!} D^\kappa\bigr)
			-(\eta+h)\bigl({}^{\gamma\!}\Dsl^\kappa\bigr)\Bigr)\;.
  \end{align*}
  Equivalently,
  \begin{align*}
	\eta_X\bigl(D^\kappa\bigr)
	&=\frac{A_G\Bigl(\frac1{\sinh\left(\frac\delta2(\punkt)\right)}
			\,e^{(\alpha-\frac\delta2)(\punkt)}\Bigr)}
		{A_G(\rho_G)}(-X)\\
	&\qquad
		+2\,S_{G/H}\Biggl(\frac{\chi_H^\kappa\bigl(e^{\punkt|_\frs}\bigr)}
				{\delta(\punkt)}
			\,\prod_{\beta\in\Delta_G^+}
				\frac{\beta(\punkt|_\frs)}{\beta(\punkt)}
			\,\prod_{\beta\in W_\pi^+}
				\frac1{2\sinh\bigl(\frac\beta2(\punkt|_\frs)
					\bigr)}\Biggr)(-X)\tag2\\
	&\qquad
		+2\int_M\Adachsl_X\bigl(TM,\nabla^0,\nabla\LC\bigr)
			\,\ch_X\bigl(V^\kappa M,\nabla^\kappa\bigr)\\
	&\qquad
		+\sum_{\gamma\in\hat G}\chi_G^\gamma\bigl(e^{-X}\bigr)
			\Bigl(\eta\bigl({}^{\gamma\!} D^\kappa\bigr)
			-(\eta+h)\bigl({}^{\gamma\!}\Dsl^\kappa\bigr)\Bigr)\;.
  \end{align*}
  If~$\rk G\ne\rk H+1$,
  then
  \begin{align*}
	\eta_X\bigl(D^\kappa\bigr)
	&=2\int_M\Adachsl_X\bigl(TM,\nabla^0,\nabla\LC\bigr)
			\,\ch_X\bigl(V^\kappa M,\nabla^\kappa\bigr)\\
	&\qquad
		+\sum_{\gamma\in\hat G}\chi_G^\gamma\bigl(e^{-X}\bigr)
			\Bigl(\eta\bigl({}^{\gamma\!} D^\kappa\bigr)
			-(\eta+h)\bigl({}^{\gamma\!}\Dsl^\kappa\bigr)\Bigr)\;.
	\tag3
  \end{align*}
  The classical $\eta$-invariant~$\eta(D^\kappa)$ is attained at~$X=0$
  in~(1)--(3).
\end{Theorem}

It is easy to see that the singularities within the parentheses of the first
term on the right hand side of~(1) cancel,
so that we are left with the alternating Weyl sum of a power series in~$X$
divided by the linearised
Weyl denominator~$\prod_{\beta\in\Delta_G^+}(-i\beta)(X)$.
The result will then be a $W_G$-invariant power series
representing a modified infinitesimal equivariant $\xi$-invariant of~$G/H$.
Note also that the last sum in~(1)--(3) is finite
because only finitely many eigenvalues
change sign when one passes from~$D^\kappa$ to~$\Dsl^\kappa$.

\begin{proof}
Assume that~$e^{-X}$ acts freely on~$M$.
Then by~\cite{G1}, \cite{G2},
	$$(\eta+h)_{e^{-X}}\bigl(D^\kappa\bigr)
		-(\eta+h)_{e^{-X}}\bigl(\Dsl^\kappa\bigr)
	=\sum_{\gamma\in\hat G}\chi_G^\gamma\bigl(e^{-X}\bigr)
		\Bigl((\eta+h)\bigl({}^{\gamma\!} D^\kappa\bigr)
		-(\eta+h)\bigl({}^{\gamma\!} \Dsl^\kappa\bigr)\Bigr)$$
is the equivariant spectral flow from~$\Dsl^\kappa$ to~$D^\kappa$.
This fits with Theorem~\ref{APSDThm},
applied to the cylinder~$M\times[0,1]$,
such that the induced operators on the ends~$M\times\{0\}$
and~$M\times\{1\}$ are precisely~$\Dsl^\kappa$ to~$D^\kappa$.
Thus by Theorem~\ref{GdreiThm},
we have
\begin{align*}
  \eta_X\bigl(D^\kappa\bigr)
  &=\eta_{e^{-X}}\bigl(D^\kappa\bigr)
	+2\int_M\frac{\thet_X}{d_X\thet_X}
		\,\Adach_X\bigl(TM,\nabla\LC\bigr)
		\,\ch_X\bigl(V^\kappa M,\nabla^0\bigr)\\
  &=(\eta+h)_{e^{-X}}\bigl(\Dsl^\kappa\bigr)
	+2\int_M\frac{\thet_X}{d_X\thet_X}
		\,\Adach_X\bigl(TM,\nabla\LC\bigr)
		\,\ch_X\bigl(V^\kappa M,\nabla^0\bigr)\\
  &\qquad
	+\sum_{\gamma\in\hat G}\chi_G^\gamma\bigl(e^{-X}\bigr)
		\Bigl(\eta\bigl({}^{\gamma\!} D^\kappa\bigr)
		-(\eta+h)\bigl({}^{\gamma\!} \Dsl^\kappa\bigr)\Bigr)\;.
\end{align*}

We still assume that~$X$ acts freely on~$M$.
Then clearly
\begin{multline*}
  \frac{\thet_X}{d_X\thet_X}
	\,\Bigl(\Adach_X\bigl(TM,\nabla\LC\bigr)
		-\Adach_X\bigl(TM,\nabla^0\bigr)\Bigr)
  =\frac{\thet_X}{d_X\thet_X}
	\,d_X\Adachsl_X\bigl(TM,\nabla^0,\nabla\LC\bigr)\\
  =\Adachsl_X\bigl(TM,\nabla^0,\nabla\LC\bigr)
	-d_X\,\biggl(\frac{\thet_X}{d_X\thet_X}
		\,\Adachsl_X\bigl(TM,\nabla^0,\nabla\LC\bigr)\biggr)\;,
\end{multline*}
so we get
\begin{align*}
	\eta_X\bigl(D^\kappa\bigr)
	&=(\eta+h)_{e^{-X}}\bigl(\Dsl^\kappa\bigr)
		+2\int_M\frac{\thet_X}{d_X\thet_X}
			\,\Adach_X\bigl(TM,\nabla^0\bigr)
			\,\ch_X\bigl(V^\kappa M,\nabla^0\bigr)\\
	&\qquad
		+2\int_M\Adachsl_X\bigl(TM,\nabla^0,\nabla\LC\bigr)
			\,\ch_X\bigl(V^\kappa M,\nabla^\kappa\bigr)\\
	&\qquad
		+\sum_{\gamma\in\hat G}\chi_G^\gamma\bigl(e^{-X}\bigr)
			\Bigl(\eta\bigl({}^{\gamma\!} D^\kappa\bigr)
			-(\eta+h)\bigl({}^{\gamma\!} \Dsl^\kappa\bigr)\Bigr)\;.
\end{align*}
Theorem~\ref{EtaThm} and Theorem~\ref{BLDThm}~(2) give~(2) if~$\rk G=\rk H+1$,
and~(3) otherwise.

To obtain~(1),
we use Theorem~\ref{BLDThm}~(1),
and we rewrite the result of Theorem~\ref{EtaThm}
using the Weyl denominator formula~\eqref{WeylFormulas}.
\end{proof}

Because the odd signature operator~$B=D$ often comes up
in topological applications,
we want to state formula~(2) for this special case.
Therefore,
let~$\pix=\kappa_1\oplus\cdots\oplus\kappa_l$
be the decomposition of~$\pix$ into $\frh$-irreducible components,
and let~$\alpha_1$, \dots, $\alpha_l$ be the corresponding weights
of~$\frg$ as in~\eqref{AlphaDef}.
Let
\begin{align*}
	\Ldach_X(TM,\nabla)
	&=\Adach_X(TM,\nabla)\wedge\ch_X(\mathcal S,\nabla)\\
	&=\det^{\textstyle\frac1{2}}
		\biggl(R^{TM}_X\,\coth\frac{R^{TM}_X}2\biggr)
\end{align*}
denote a rescaled version of Hirzebruch's $L$-genus,
and let~$\Ldachsl_X$ denote the corresponding equivariant Chern-Simons
class.

\begin{Corollary}\label{BCor}
  Assume that~$\rk G=\rk H+1$.
  Then the infinitesimally equivariant $\eta$-invariant
  of the odd signature operator~$B$ on~$M$ is given by
  \begin{align*}
    \eta_X(B)
    &=\sum_l\frac{A_G\Bigl(\frac1{\sinh\left(\frac\delta2(\punkt)\right)}
		\,e^{(\alpha_l-\frac\delta2)(\punkt)}\Bigr)}
	{A_G(\rho_G)}(-X)\\
    &\qquad
	+2\,S_{G/H}\Biggl(
		\frac1{\delta(\punkt)}
		\,\prod_{\beta\in\Delta_G^+}
			\frac{\beta(\punkt|_\frs)}{\beta(\punkt)}
		\,\prod_{\beta\in W_\pi^+}
			\coth\biggl(\frac{\beta(\punkt|_\frs)}2\biggr)
	\Biggr)(-X)\tag1\\
    &\qquad
	+2\int_M\Ldachsl_X\bigl(TM,\nabla^0,\nabla\LC\bigr)
	+\sum_{\gamma\in\hat G}\chi_G^\gamma\bigl(e^{-X}\bigr)
		\Bigl(\eta\bigl({}^{\gamma\!} B\bigr)
		-(\eta+h)\bigl({}^{\gamma\!}\tilde B\bigr)\Bigr)\;.
  \end{align*}
  If~$\rk G\ne\rk H+1$, then
  \begin{align*}
    \eta_X(B)
    &=2\int_M\Ldachsl_X\bigl(TM,\nabla^0,\nabla\LC\bigr)
	+\sum_{\gamma\in\hat G}\chi_G^\gamma\bigl(e^{-X}\bigr)
		\Bigl(\eta\bigl({}^{\gamma\!} B\bigr)
		-(\eta+h)\bigl({}^{\gamma\!}\tilde B\bigr)\Bigr)\;.\tag2
  \end{align*}
  The classical $\eta$-invariant~$\eta(B)$ is attained at~$X=0$
  in both cases.
\end{Corollary}

\begin{proof}
  Let~$\beta_1$, \dots, $\beta_r$ denote the positive weights of~$\pi$,
  counted with the right multiplicity.
  If~$\rk G=\rk H+1$,
  then the weights of~$\pix$ take
  the form~$\bigl(\pm\frac{\beta_1}2\pm\dots\pm\frac{\beta_r}2\bigr)|_\frs$,
  each with multiplicity one~\cite{G1}, \cite{G2}.
  Therefore,
	$$\chi_H^\pix\bigl(e^{-X|_\frs}\bigr)
	=\prod_{\beta\in W_\beta^+}
		2\cosh\biggl(\frac\beta2(-X|_\frs)\biggr)\;.$$

  Note also that for the odd signature operator~$B$,
  we need the Levi-Civita connection on the twist bundle
  instead of the reductive connection.
  This is why the correct Chern-Simons contribution
  is given by~$\Ldachsl_X(TM,\nabla^0,\nabla\LC)$,
  not by~$\Adachsl_X(TM,\nabla^0,\nabla\LC)\,\ch_X(\mathcal S,\nabla^0)$.
  The Corollary now follows easily from Theorem~\ref{MainTheorem}~(2) and~(3).
\end{proof}

\bibliographystyle{alpha}

\enddocument